\documentclass[10pt]{amsart}
\usepackage{latexsym,amsmath,amscd, amssymb}
\allowdisplaybreaks \headheight=8pt     \topmargin=0pt
\newtheorem{definition}{Definition}

\newtheorem{proposition}{Proposition}
\newtheorem{corollary}{Corollary}
\newtheorem{lemma}{Lemma}
\theoremstyle{remark}

 \title{D\MakeLowercase{rinfel'd}-I\MakeLowercase{hara}
R\MakeLowercase{elations for the} C\MakeLowercase{rystalline} F\MakeLowercase{robenius} }
\date{September 10, 2003}
\author{\bf S\MakeLowercase{inan} \"{U}\MakeLowercase{nver}}
\begin{document} \maketitle \noindent

{\it Notation and convention.} If we have an element $a$ of the
ring of associative formal power series over a ring $A,$ in the
variables $x_i,$ $i \in I,$ we denote the coefficient of $x^{J}$
in $a,$ for a multi-index $J,$ by $a[ x^{J} ].$ By a variety $X$
over a field $K,$ we mean a geometrically integral  $K$-scheme
$X,$ that is separated and of finite type over $K.$

\section{Introduction}

For $K$ a field of characteristic zero,  let ${\rm M}(K)$ denote
the set of Drinfel'd associators defined over $K$ [Dr]. The
variety $M$ is of great interest because of its connection to the
deformations of universal enveloping algebras, fundamental group
of the Teichmuller tower, and  to ${\rm
Gal}(\overline{\mathbb{Q}}/\mathbb{Q}).$ There is a natural map
${\rm M}\to \mathbb{A}^1$ and for $\lambda \in K,$ and ${\rm
M}_{\lambda}$ is a torsor under ${\rm GRT}_{1}\, (={\rm M}_{0}).$
If $K<< X,Y >>$ denotes the formal associative power series in $X$
and $Y$ over $K$  then ${\rm GRT}_{1}(K)$ is defined to be the set
of elements $\varphi  \in K<<X,Y >>$ that satisfy

\begin{eqnarray*}
&(i)&\;\; \varphi(X,Y)\cdot \varphi(X,Y)=1\\
&(ii)&\;\; \varphi(Z,X)\cdot \varphi(Y,Z) \cdot \varphi(X,Y)=1,
\;\;{\rm when  \;\;} X+Y+Z=0\\
&(iii)& \varphi(x_{23},x_{34})\cdot \varphi(x_{40},x_{01})\cdot
\varphi(x_{12},x_{23}) \cdot \varphi(x_{34},x_{40})\cdot
\varphi(x_{01},x_{12})=1,\\
& &{\rm \; when \;} x_{ii}=0,\; \sum_{j}x_{ij}=0, \; {\rm and\;}
[x_{ij},x_{kl}]=0, \; {\rm if \;} \{i,j \}\cap \{k,l \}=\emptyset.
\end{eqnarray*}

The multi-zeta values (or equivalently the K-Z equation) gives an
element in ${\rm GRT}_1(\mathbb{C})$ [Dr], the Galois action on
$\pi_{1,et}(\mathbb{P}^{1}\setminus \{0,1,\infty \},\cdot)$ gives
an element in ${\rm GRT}_{1}(\mathbb{Q}_{\ell})$ [Ih]. The aim  of
the following is to show that the crystalline frobenius on
$\pi_{1,dR}(\mathbb{P}^1 \setminus \{0,1,\infty \},\cdot )$ also
gives an element ${\rm GRT}_{1}(\mathbb{Q}_p).$ The main part is
the proof of (iii) which is an application of the theory of
tangential basepoints. Therefore most of the following is devoted
to developing the crystalline  theory of tangential basepoints.
The Betti, de Rham, and \'{e}tale realizations of the tangential
basepoints in the case of curves  were given in the fundamental
paper of Deligne [De].

{\it Acknowledgements.} The author would like to thank Professor
Arthur Ogus for his  suggestions on the tangential basepoints and for his encouragement.

\section{ Pullback to the log point} Let $\overline{X}/k$ be a
smooth variety over a field $k.$ Let $D \subseteq \overline{X}$ be
a simple normal crossings divisor, $X:=\overline{X}\setminus D$
and $x \in D(k).$ Let $\{ D_{i}\}_{i \in I},$ be the set of
irreducible components of $D$ passing through $x,$ $D_x:=\cup_{i
\in I }D_i,$ $I_{i}$ the ideal of $D_{i}$ in $\mathcal{O}_{X,x},$
and $d: \mathfrak{m}_{x}\to \mathfrak{m}_x/\mathfrak{m}_x^2$ the
canonical projection. Note that $I_i \subseteq \mathfrak{m}_{x}.$
Denote by $\overline{X}_{log}$ the canonical fine saturated log
scheme (in the Zariski topology) defined by $(\overline{X},D).$
The underlying scheme of $\overline{X}_{log}$ is $\overline{X}$
and the log structure is defined by the inclusion
$$
M_{\overline{X}}:=\mathcal{O}_{\overline{X}} \cap
i_{*}(\mathcal{O}_{X}^*) \to \mathcal{O}_{\overline{X}},
$$
where $i:X \hookrightarrow \overline{X}$ is the inclusion. Note
that
$$\overline{M}_{\overline{X},x}:=
M_{\overline{X},x}/\mathcal{O}^*_{\overline{X},x}\simeq {\rm
Cart}^{-}(\overline{X},D_x),
$$
canonically, where  ${\rm Cart}^{-}(\overline{X},D_x)$ denotes the
monoid of anti-effective Cartier divisors on $\overline{X}$
supported on $D_x.$ If the $D_i$ are  defined locally by $t_i=0,$
with $t_i \in \mathcal{O}_{\overline{X},x},$ then
$M_{\overline{X},x}=(\mathcal{O}_{\overline{X},x})^* _{\prod
t_i}.$

 Let $x_{log}$ denote the log scheme obtained by
pulling back the log structure on $\overline{X}_{log}$ via the map
$ {\rm Spec}k \to \overline{X} $ corresponding to $x.$ Note that
the monoid $M_{x_{log}}$ on $x_{log}$ is $M_{\overline{X},x}
\otimes _{\mathcal{O}^*_{\overline{X},x}}k^*.$

{\it Notation.} With the notation above let
$$C(D_x,\overline{X}):=\{\varphi: \varphi {\rm \; is\; a\; splitting\;
of\; } M_{x_{log}} \to {\rm Cart}^{-}(\overline{X},D_{x}) \}.
$$
Then
$$C(D_x,\overline{X})\simeq
\{ (\cdots, \overline{t}_i,\cdots): \overline{t}_i \in d I_i
\setminus \{0 \}, {\rm \;\; for\;\;} i \in I\}.$$ We say that
$(\cdots,\overline{t}_i,\cdots)$  is {\it transversal} to $D_x$ if
$(\cdots,\overline{t}_i,\cdots) \in C(D_x,\overline{X}).$
Similarly let
$$
N_{D_x/\overline{X}}:=\prod _{i \in I} N_{D_{i}/\overline{X}}(x),
$$ be
the fiber at $x$ of the product of the normal bundles of $D_i$ in
$\overline{X}.$ Note that elements of $C(D_x,\overline{X})$ are in
one to one correspondence with linear isomorphisms
$$
\prod _{i \in I}N_{D_i/\overline{X}}(x) \simeq \prod _{i \in I
}\mathbb{A}^1
$$
preserving the factors, and hence with
$N_{D_x/\overline{X}}^*(k),$ where we let
$N_{D_x/\overline{X}}^*:=\prod_{i \in I}N_{D_i/\overline{X}}^*,$
with $N_{D_{i}/\overline{X}}^*:=N_{D_{i}/\overline{X}}\setminus \{
0\}.$ We call elements $(\cdots, v_i,\cdots )\in
N^*_{D_x/\overline{X}}(k),$ r-tuples of vectors transversal to
$D_x.$ Therefore

\begin{lemma}
The set of splittings of the log structure on $x_{log}$ is in one
to one correspondence with  $N^*_{D_x/\overline{X}}(k).$
\end{lemma}

Let ${\rm Spec}\, k_{x,log}$ denote the log scheme with underlying
scheme ${\rm Spec} \, k, $ and log structure associated to the
pre-log structure ${\rm Cart}^{-}(\overline{X},D_x) \to k$ that
maps all the nonzero elements of ${\rm
Cart}^{-}(\overline{X},D_x)$ to 0. From the lemma above the set of
isomorphisms $ x_{log} \simeq {\rm Spec}k_{x,log}$ is in one to
one correspondence with $N_{D_{x}/\overline{X}}^*(k).$

\section{ Construction of tangential basepoints in the logarithmic
case}

{\bf 3.1. Crystalline.} Let $k$ be a perfect field of
characteristic p.  Let $W:=W(k)$ be the ring of Witt vectors over
$k,$ with field of fractions $K,$  and let $S:={\rm Spec}W.$
Associated to the exact closed immersion
$$
x_{log} \to \overline{X}_{log}
$$
we have a pullback functor
$$
x_{log}^*:{\rm Isoc}(\overline{X}_{log}/S) \to {\rm
Isoc}(x_{log}/S)
 $$
between the categories of convergent log isocrystals on
$\overline{X}_{log}$ and $x_{log}.$ Choosing an element $v$ in
$N_{D/\overline{X}}^*(k)$ gives a functor
$$
v^*:{\rm Isoc}(\overline{X}_{log}/S) \to {\rm Isoc}({\rm Spec}\, k
_{x,log}/S).
 $$
\begin{definition} Let $\mathcal{T}_x/K$ denote the tannakian category of
vector bundles $V$ over $K$ endowed with a homomorphism
$\varphi:{\rm Cart}^{-}(\overline{X},D_x)^{gp\, *} \to {\rm
End}(V),$ where  $${\rm Cart}^{-}(\overline{X},D_x)^{gp\, *}:={\rm
Hom }_{\mathbb{N}}({\rm Cart}^-(\overline{X},D_x ),\mathbb{Z}).$$
The tensor product of $(V,\varphi)$ and $(W,\psi)$ is defined to
be $V \otimes W$ endowed with the map that sends $-D_{i} \in {\rm
Weil}^{-}(\overline{X},D_{x} )$ to $id \otimes \psi(-D_i)+
\varphi(-D_i)\otimes id \in {\rm End}(V \otimes W).$
 \end{definition}
There is a natural realization functor
$$
ev:{\rm Isoc}({\rm Spec}\, k_{x,log}/S) \to \mathcal{T}_x
$$
defined as follows. Let $S_{x,log}$ denote ${\rm Spec}W[{\rm
Cart}^-(\overline{X},D_x)]$ with the standard log structure. Note
that $S_{x,log}/S$ is smooth and there is a canonical exact closed
immersion
$$
{\rm Spec} k_{x,log} \to S_{x,log}.
$$
If $(E,\nabla) \in {\rm Isoc}({\rm Spec}k_{x,log}/S)$ then we have
the realization $$(E,\nabla)_{S_{x,log}} \in {\rm
Mic}(S_{x,log,K}/K).$$ The connection induces a map
$$
-R: E(0) \to E(0) \otimes \Omega^{1}_{S_{x,log,K}/K}(0).
$$
Since $ \Omega^{1}_{S_{x,log,K}/K}(0) \simeq {\rm
Cart}^{-}(\overline{X},D_x)^{gp} \otimes _{\mathbb{Z}} K,$ we
obtain an object $$((E,\nabla)_{S_{x,log, K}}(0),R ) \in
\mathcal{T}_x/K$$ This defines the functor above. We let $ev(v):=
 ev  \circ v^*,$ and $ev(v)$  composed with the natural fiber
 functor of $\mathcal{T}_x$ to be
$$
fib(v): {\rm Isoc} (\overline{X}_{log}/W) \to Vec_{K}
$$
over $K,$ for each  $v \in N_{D/\overline{X}}^*(k), $ where
$Vec_{K}$ denotes the category of finite dimensional vector spaces
over $K.$

{\bf 3.2. de Rham.} Let $k$ be a field of characteristic zero.
Then we have a canonical functor, which does not depend on the
choice of a set of  vectors transversal to $D,$
$$
ev_{dR}:{\rm  Mic}(x_{log}/k) \to \mathcal{T}_x/k
$$
defined as follows. Let  $(E,\nabla)\in {\rm Mic}(x_{log}/k).$
Then as above we have a residue map
$$
R:E \to E\otimes \Omega^{1}_{x_{log}/k},
$$
and an isomorphism $\Omega^{1}_{x_{log}/k} \simeq {\rm
Cart}^{-}(\overline{X},D_{x})^{gp} \otimes_{\mathbb{Z}}k$ which
gives an object of $\mathcal{T}_{x}$ as above. For $v \in
N_{D_x/\overline{X}} ^* (k) $ we obtain a commutative diagram
\[
\begin{CD}
{\rm Mic}(x_{log}/k) @>ev_{dR}>> \mathcal{T}_x\\
 @VV{id_{v}}V @VV{id}V \\
{\rm Mic}(k_{x,log}/k) @>ev_{dR}>> \mathcal{T}_x,
\end{CD}
\]
where $id_{v}$ denotes the identification as above that depends on
the choice of a transversal set of tangent vectors; and where the
lower horizontal map is defined analogously to the horizontal map
above. Let $ev_{dR}(x):= ev_{dR}  \circ x_{log} ^*$ and
$ev_{dR}(v):=   ev_{dR} \circ id_{v} \circ x_{log}^*;$ and
$fib_{dR}(\cdot),$ $e_{dR}(\cdot)$ composed with the standard
fiber functor of $Vec_{k}.$ Then by the above we have a canonical
isomorphism $fib_{dR}(x) \simeq fib_{dR}(v)$ and in particular, in
the logarithmic de Rham theory the fiber functor is independent of
the choice of vectors.

{\bf 3.3. Comparison.} With the notation as in the beginning of
this section, assume that $k$ is a perfect field of characteristic
$p.$ Fix $v \in N_{D_x/\overline{X}}^*(k),$ a vector   transversal
to $D_x.$ Let $\overline{\mathfrak{X}}/W$ be a smooth formal
scheme, $\mathfrak{D} \subseteq \overline{\mathfrak{X}}$ a simple
relative normal crossings divisor, $\mathfrak{x} \in
\mathfrak{D}(W)$ such that the reductions of
$\overline{\mathfrak{X}},$ $\mathfrak{D},$ and $\mathfrak{x}$ are
$\overline{X},$ $D,$ and $x$ respectively. We let $\mathfrak{X} :=
\overline{\mathfrak{X}}\setminus \mathfrak{D},$
$\overline{\mathfrak{X}}_{log}$ the log scheme associated to the
pair $(\overline{\mathfrak{X}},\mathfrak{D})$ etc.  Then
$N_{\mathfrak{D}_{\mathfrak{x}}/\overline{\mathfrak{X}}}$ is a
lifting of $N_{D_x/\overline{X}}.$ Let $\mathfrak{v}\in
N_{\mathfrak{D}_{\mathfrak{x}} /\overline{\mathfrak{X}}}^*(W),$
with reduction $v.$   We have a commutative diagram
\[
\begin{CD}
{\rm Spec}\, k_{x,log} @>>> {\rm Spec}\, W_{\mathfrak{x},log}\\
@VV{id_{v}}V  @VV{id_{\mathfrak{v}}}V\\
{\rm Spec}\, x_{log} @>>> {\rm Spec}\, \mathfrak{x}_{log}\\
@VVV @VVV\\
\overline{X}_{log} @>>> \overline{\mathfrak{X}}_{log},\\
\end{CD}
\]
and a canonical exact closed imbedding ${\rm Spec}
W_{\mathfrak{x},log} \hookrightarrow
\mathfrak{S}_{\mathfrak{x},log}.$ Note that ${\rm Spec}
W_{\mathfrak{x},log}$ does not depend on $\mathfrak{x},$ therefore
we will denote it by ${\rm Spec}W_{x,log}.$

For $(E,\nabla) \in {\rm Isoc}(\overline{X}_{log}/W),$ let
$(E,\nabla)_{\overline{\mathfrak{X}}} \in {\rm
Mic}(\overline{\mathfrak{X}}_{K,log}/K )$ denote the realization
of $(E,\nabla)$ on the log rigid analytic space
$\overline{\mathfrak{X}}_{K,log}.$ Then since $id_{\mathfrak{v}}$
is the identity map on the underlying schemes, we have
$fib(v)(E,\nabla)\simeq
E_{\overline{\mathfrak{X}}}(\mathfrak{x}_K)$ by the definition of
$fib(v).$We call this the {\it realization} of the fiber functor
$fib(v)$ corresponding to the data
$(\overline{\mathfrak{X}},\mathfrak{D},\mathfrak{x},\mathfrak{v}
),$ and denote it by $fib(v)_{\mathfrak{v}}$ in order to remember
the choice of the model. Note that this in fact does not depend on
the choice of the lifting $\mathfrak{v}$ as expected from the de
Rham version. However the comparison with the different models,
i.e. the crystalline version, depends on the choice of the
liftings of the tangent vectors. Fixing the model gives the
comparison between the de Rham and the crystalline versions.

Let $(\overline{\mathfrak{Y}},\mathfrak{E},
\mathfrak{y},\mathfrak{u})$ be another data of a lifting. Let
$(\overline{\mathfrak{X}} \times \overline{\mathfrak{Y}} )^{\sim}$
denote the blow-up of $\overline{\mathfrak{X}} \times
\overline{\mathfrak{Y}}$ along  $\cup_i (\mathfrak{D}_i \times
\mathfrak{E}_i).$ Denote $\mathfrak{D}_i \times \mathfrak{E}_i$ by
$\mathfrak{F}_i.$ Endow $(\overline{\mathfrak{X}}\times
\overline{\mathfrak{Y}})^{\sim}$ with the log structure associated
to the exceptional divisor of the blow-up, and denote it by
$(\overline{\mathfrak{X}}\times
\overline{\mathfrak{Y}})_{log}^{\sim}.$ The fiber of the blow-up
over the point $(\mathfrak{x}, \mathfrak{y})$ is isomorphic to
$$
\times_i
\mathbb{P}(N_{\mathfrak{F}_{i}/\overline{\mathfrak{X}}\times
\overline{\mathfrak{Y}} } ),
$$
the product of the projective normal bundles of $\mathfrak{F}_i$
in $\overline{\mathfrak{X}} \times \overline{\mathfrak{Y}}.$ In
particular the pair $[\mathfrak{v},\mathfrak{u}]$ defines an
element of the fiber of
$$
(\overline{\mathfrak{X}} \times \overline{\mathfrak{Y}} )^{\sim}
\to \overline{\mathfrak{X}} \times \overline{\mathfrak{Y}}
$$
over $(\mathfrak{x},\mathfrak{y}),$ where if
$\mathfrak{v}=(\cdots,\mathfrak{v}_{i},\cdots ),$ and
$\mathfrak{u}=(\cdots,\mathfrak{u}_i,\cdots)$ then
$$[\mathfrak{v},\mathfrak{u}]
:=(\cdots, [\mathfrak{v}_i,\mathfrak{u}_i],\cdots).$$

If $(E,\nabla) \in {\rm Isoc }(\overline{X}_{log}/W)$ then we have
a canonical isomorphism between the pull-backs of
$(E,\nabla)_{\overline{\mathfrak{X}}}$ and
$(E,\nabla)_{\overline{\mathfrak{Y}}}$ to
$(\overline{\mathfrak{X}} \times
\overline{\mathfrak{Y}})_{K}^{\sim},$ and hence evaluating this
isomorphism at $[\mathfrak{v}, \mathfrak{u} ]_K$ gives an
isomorphism
$$
E_{\overline{\mathfrak{X}}}(\mathfrak{x}_K)\simeq
E_{\overline{\mathfrak{Y}} }(\mathfrak{y}_K).
 $$
On the other hand from the above  diagram applied to
$\overline{\mathfrak{X}}, \mathfrak{x}, \mathfrak{v}$ and
$\overline{\mathfrak{Y}}, \mathfrak{y}, \mathfrak{u},$ we obtain
isomorphisms, functorial in $E,$
$$
E_{\overline{\mathfrak{X}}}(\mathfrak{x}_K) \simeq
fib(v)_{\mathfrak{v}}(E,\nabla) \simeq
fib(v)_{\mathfrak{u}}(E,\nabla)\simeq
E_{\overline{\mathfrak{Y}}}(\mathfrak{y}_K).
$$

{\bf Claim.} The two isomorphisms above, between
$E_{\overline{\mathfrak{X}}}(\mathfrak{x}_K)$ and
$E_{\overline{\mathfrak{Y}}}(\mathfrak{y}_{K}),$  are the same.

{\it Proof.} Let $M:=M_{W_{x,log}},$ and $\Delta: M \oplus M \to
M,$ the map induced by $(id_{M},id_M),$ and $I$ the ideal $
\Delta^*((-\sum_{i \in I} D_i)).$ Let $(W_{x,log}\times
W_{x,log})^{\sim}$ denote the blow-up of $W_{x,log}\times
W_{x,log}$ along $I.$ Note that its underlying scheme is $\prod_{i
\in I }\mathbb{P}_{W}^1.$

Let
$$\alpha_{\mathfrak{v}}:{\rm Spec}W_{x,log} \to
\overline{\mathfrak{X}}$$ and
$$\alpha_{\mathfrak{u}}:{\rm Spec}W_{x,log} \to
\overline{\mathfrak{Y}}$$ be the exact closed imbeddings induced
by $\mathfrak{v}$ and $\mathfrak{u}.$ These induce a map
$$
(W_{x,log}\times W_{x,log})^{\sim} \to (\overline{\mathfrak{X}}
\times \overline{\mathfrak{Y}} )_{log}^{\sim}
$$
that maps $(1,\cdots,1)$ to $[\mathfrak{v},\mathfrak{u}].$ The
diagonal map factors through the blow-up to give a map $\Delta:
{\rm Spec} W_{x,log} \to (W_{x,log}\times W_{x,log})^{\sim},$
which has image $(1,\cdots,1).$ This map makes the diagrams

\[
\begin{CD}
{\rm Spec} W_{x,log} @>>> {\rm Spec} (W_{x,log} \times
W_{x,log})^{\sim}@>>> (\overline{\mathfrak{X}}\times
\overline{\mathfrak{Y}})_{log}^{\sim}\\
@VV{id}V  & &   @VVV\\
 {\rm Spec} W_{x,log}  &  @>{\alpha_{\mathfrak{v}}}>> &
\overline{\mathfrak{X}}_{log}
\end{CD}
\]
and

\[
\begin{CD}
{\rm Spec} W_{x,log} @>>> {\rm Spec} (W_{x,log} \times
W_{x,log})^{\sim}@>>> (\overline{\mathfrak{X}}\times
\overline{\mathfrak{Y}})_{log}^{\sim}\\
@VV{id}V  & &   @VVV\\
 {\rm Spec} W_{x,log}  &  @>{\alpha_{\mathfrak{u}}}>> &
 \overline{\mathfrak{Y}}_{log}
\end{CD}
\]
commute. Therefore if $(E,\nabla) \in {\rm
Isoc}(\overline{X}_{log}/W)$ then the isomorphism between
$$fib(v)_{\mathfrak{v}}(E,\nabla) {\rm \;\;\; and
\;\;\; } fib(v)_{\mathfrak{u}}(E,\nabla)$$  is obtained by pulling
back the canonical isomorphism between $$p_1^*
(E,\nabla)_{\overline{\mathfrak{X}}} {\rm \;\; and \;\;}
p_2^*(E,\nabla)_{\overline{\mathfrak{Y}}}$$ via
$$
{\rm Spec}W_{x,log} \to (\overline{\mathfrak{X}}\times
\overline{\mathfrak{Y}} )^{\sim}_{log},
$$
where the $p_i,$ $i=1,2$ denote the projections. Since the last
map has image $[\mathfrak{v},\mathfrak{u}],$ this proves the
claim. \hfill $\Box$

 \section{  Comparison of different log points} Let
 $D_{(x)} \subseteq N_{D_x/\overline{X}}$ denote the simple normal
 crossings divisor defined by the ideal $\prod_{i \in I}dI_i.$
In other words if $t_1, \cdots ,t_n$ is a regular system of
parameters on $\overline{X}$ at $x,$ such that $D$ is defined
locally by $t_1 \cdots t_r =0$ then $D_{(x)}$ is defined by
$d{t}_1 \cdots d{t}_{r}=0.$ Note that $\{d{t}_i\}_{1\leq i \leq
r}$ is a system of coordinates for $N_{D_{x}/\overline{X}}.$ Let
$N_{log}:=N_{D_{x}/\overline{X},log}$ be the scheme
$N_{D_{x}/\overline{X}}$ endowed with the logarithmic structure
associated to the divisor $D_{(x)},$ and $0_{log}$ the logarithmic
scheme induced by the inclusion $0 \to N_{D_{x}/\overline{X}}$
which gives an   exact closed immersion $$ 0_{log} \to N_{log}.
$$

\begin{lemma}
The log schemes $x_{log}$ and $0_{log}$ are canonically
isomorphic.
\end{lemma}

{\it Proof.} Note that both log schemes have the same underlying
scheme ${\rm Spec}\, k.$

Let $\{ t_{i}\}_{i \in I}$ be a part of a system of parameters for
$\mathfrak{m}_x$ such that  $D$ is defined locally by $t_1 \cdots
t_r=0.$ This gives a map
$$
\varphi: \otimes _{1\leq i \leq r}
 Sym^{\cdot} (dI_{i}/d I_{i}^{2}) \to \mathcal{O}_{\overline{X},x}, $$
with the property $\varphi(dt_i)=t_i.$ Since
$M_{N}=(\mathcal{O}_{N,0})_{d{t}_1 \cdots d{t}_r} ^*,$
$M_{\overline{X},x}=(\mathcal{O}_{X,x})_{t_1 \cdots t_r}^*,$ this
gives a map $M_{N} \to M_{\overline{X},x}$ and hence a map
$$
M_{T}\otimes_{\mathcal{O}_{N,x}^*}k^* \to
M_{\overline{X},x}\otimes_{\mathcal{O}_{\overline{X},x}^*} k^*
$$
that we continue to denote by $\varphi.$

The last map is independent of the choice of $\{t_{i} \}_{1\leq
i\leq r}$ as above. If $\{s_i \}_{i \in I}$ is another set then
$s_i/t_i \in \mathcal{O}_{\overline{X},x}^*,$ for $i \in I.$ Let
$\psi$ be the map corresponding to $\{ s_i\}_{i \in I}.$ To show
the independence it is enough to show that $\psi$ and $\varphi$
agree on $\{ds_i \}_{i \in I}.$ Note that
$$
\varphi(ds_i)=\varphi(d (\frac{s_i}{t_i} \cdot t_i)
)=\varphi(d(\frac{s_i}{t_i})\cdot t_i + \frac{s_i}{t_i}(x) \cdot
dt_i )=\frac{s_i}{t_i}(x)\cdot t_i
$$
and $\psi(ds_i)=s_i.$ Therefore
$$\frac{\varphi(ds_i)}{\psi(ds_i)}=\frac{s_i}{t_i}(x) \cdot
\frac{t_i}{s_i}$$ and hence $\varphi(ds_i)$ and $\psi(ds_i)$ are
equal in
$M_{\overline{X},x}\otimes_{\mathcal{O}_{\overline{X},x}^*}k^*, $
for $i \in I.$ This shows the independence.

This map makes the diagram
\[
\begin{CD}
0 @>>> k^* @>>> M_{0_{log}} @>>> \overline{M}_{0_{log}} @>>> 0\\
& & @VV{id}V @VVV @VV{\alpha_x}V\\
0 @>>> k^* @>>> M_{x_{log}} @>>> \overline{M}_{x_{log}}@>>> 0
\end{CD}
\]
commute, where $\alpha_x$ is the canonical isomorphism
$$\alpha_{x}:{\rm Cart}^{-}(N_{D_{x}/\overline{X}},D_{(x)})
\simeq {\rm Weil}^-(N_{D_{x}/\overline{X}},D_{(x)}) \to {\rm
Weil}^-(\overline{X},D_{x})\simeq {\rm
Cart}^{-}(\overline{X},D_{x}).$$ This implies that the middle map
is an isomorphism. \hfill $\Box$

\section{ Tangential basepoints in the  regular singular case} We
continue with the notation above. Let $ \hat{N }_{log}$ and
$\hat{\overline{X}}_{log}$ denote the completions of $N_{log}$
along 0 and $\overline{X}$ along $x.$ Let $\psi:
  \hat{N}_{log} \to \hat{\overline{X}}_{log}  $ be a morphism making
the diagram
\[
\begin{CD}
 0_{log} @>>>         x_{log} \\
@VVV @VVV \\
\hat{N}_{log}    @>\psi>>    \hat{\overline{X}}_{log}
\end{CD}
\]
commute. Such a $\psi$ can be obtained by choosing  a system of
parameters $t_1, \cdots, t_n$ at $x$ on $\overline{X}$ such that
$D_{x}$ is locally defined by $t_1 \cdots t_r=0,$ and letting
$\psi^*(t_i)=dt_i \cdot\varphi_i$ with $\varphi_i \in 1 +
\hat{\mathfrak{m}}_{0} ,$ for $1 \leq i\leq r,$ and
$\psi^*(t_{j})=s_j,$
 with  $s_j \in \hat{\mathfrak{m}}_{0} ,$ for $r+1 \leq j \leq
 n.$

{\bf 5.1. Construction.} In this section assume that the base
field $K$ is of characteristic zero.  We say that an integrable
connection $\nabla$ on a vector bundle $E$ on $K[[t_1,\cdots,t_n
]]_{t_1\cdots t_r}$ is regular along the divisor $t_1 \cdots
t_r=0,$ if $(E,\nabla)$ has a logarithmic  extension
$(\overline{E},\nabla)$ to $K[[t_1,\cdots, t_n]].$ If $Y/K$ is
smooth, and $D \subseteq Y$ is a simple normal crossings divisor
then we denote the category of vector bundles with integrable
connection on $Y \setminus D$ with regular singularities along $D$
by  ${\rm Mic}_{ reg, D}(Y/K).$

We define a natural functor,
$$
\hat{\varphi}^*: {\rm Mic}_{reg,\hat{D}_x}(\hat{\overline{X}}/K)
\to {\rm Mic}_{reg,\hat{D}_{(x)}}(\hat{N}/K).
$$
We would like to define $\hat{\varphi}^*$ as $\psi^*.$ However
since there are different choices for $\psi,$ we have to give
canonical isomorphisms between $\psi^*$ and $\psi '\, ^*$ for two
such  choices $\psi$ and $\psi',$ satisfying the cocycle condition
for three choices. This is done exactly as above.  Let
$(\hat{\overline{X}}\times\hat{\overline{X}})^{\sim}$ be the
blow-up of  $(\hat{\overline{X}}\times\hat{\overline{X}})$ along
$\cup_{i \in I }(\hat{D}_{i} \times \hat{D}_{i}).$ Let
$\hat{D}_{x} ^{\sim}$ denote the exceptional divisor. Then $\psi$
and $\psi '$ induce a map $(\psi,\psi')^{\sim}: \hat{N}_{log} \to
(\hat{\overline{X}} \times \hat{\overline{X}} )_{log} ^{\sim}$
that factors through the completion
$\hat{\Delta}_{\hat{\overline{X}}} ^{\sim}$ of
$(\hat{\overline{X}} \times \hat{\overline{X}} )_{log} ^{\sim}$
along the strict transform $\Delta_{\hat{\overline{X}}} ^{\sim}$
of the diagonal.

Let $(E,\nabla) \in {\rm
Mic}_{reg,\hat{D}_x}(\hat{\overline{X}}/K)$ then it has  a
logarithmic extension $(\overline{E},\nabla).$ The connection
gives a horizontal isomorphism between $p_{1}
^{*}(\overline{E},\nabla)$ and $p_2 ^{*} (\overline{E},\nabla)$ on
the first infinitesimal neighborhood of the strict transform of
the diagonal in $(\hat{\overline{X}} \times \hat{\overline{X}} )
^{\sim}.$ Since the connection is integrable and $K$ has
characteristic 0 the isomorphism extends to the formal
neighborhood. The restriction of this isomorphism to
$(\hat{\overline{X} } \times \hat{\overline{X}})^{\sim}\setminus
\hat{D}_{x} ^{\sim}$ depends only on $(E,\nabla).$ Therefore
pulling back  this isomorphism via
$$
\hat{N} \setminus \hat{D}_{(x)} \to (\hat{\overline{X} } \times
\hat{\overline{X}})^{\sim}\setminus \hat{D}_{x} ^{\sim}
$$
gives the isomorphisms between $\psi^*$ and $\psi ' \, ^*$ we were
looking for. As usual the integrability of the connection gives
the cocycle condition.

We will need the following.

\begin{lemma}
If $(E,\nabla) \in {\rm Mic}_{reg,t_1 \cdots t_r }(K[[t_1,\cdots,
t_n  ]]/K)$ then it has   an extension of the form
$(\mathcal{O}^{\oplus rk E},d- \sum_{1\leq  i \leq r}
\Gamma_{i}\frac{dz_i}{z_i} ) ,$  to $K[[t_1,\cdots, t_ n ]],$
where $\Gamma_{i} \in {\rm M}_{rk E \times rkE}(K) $ and the
$\Gamma_i$ do not have eigenvalues that differ by a non-zero
integer.
\end{lemma}

{\it Proof.} This is exactly Th\'{e}or\`{e}me  3.4. in [GL] except
that there $K$ is assumed to be algebraically closed. The same
proof works without this assumption as follows. Using the notation
in loc. cit. let $\Lambda$ be a lattice, not necessarily free, in
$\Gamma(K[[t_1,\cdots,t_{n} ]]_{t_1 \cdots t_r},E )$ relative to
which $\nabla$ has logarithmic singularities along $t_1 \cdots
t_r=0.$ Let $W$ be a $K$-subspace of $\Lambda$ invariant under the
semi-simplifications $^{s}\nabla_{t_i
\partial/ \partial t_i}$ of the  $\nabla_{t_i \partial/ \partial
t_i}$ and complementary to $\mathfrak{m} \Lambda.$ Since $^s
\nabla_{t_i \partial/ \partial t_i}$ are commuting semi-simple
operators, $W$ has a decomposition $W=\oplus_{j\in J} W_j,$ as a
$K[^s \nabla_{t_i \partial/
\partial t_i}]_{1 \leq i \leq n}$ module, such that each $^s \nabla_{t_i
\partial/
\partial t_i}|_{W_j}$ is annihilated by an irreducible polynomial.
Assume that, $^s \nabla_{t_i \partial/
\partial t_i},$ for some $i,$ has two eigenvalues $\lambda_{1}$ and
$\lambda_{2}$ such $\lambda_{1}-\lambda_{2}=n \in \mathbb{Z}^{+}.$
Note that these eigenvalues cannot be the root of the same
irreducible polynomial $p(x),$ since otherwise $\lambda_1$ will be
a root of $p(x)-p(x-n),$ a non-zero polynomial of degree less than
$p(x).$ If $\lambda_{k},$ $1\leq k\leq 2,$ appear as an eigenvalue
of $^s \nabla_{t_i \partial/ \partial t_i}$ on $W_{j_k},$ then the
corresponding minimal  polynomials satisfy
$p_{i_1}(x)=p_{i_2}(x-n).$ We replace $W_{i_2}$ with $t_{i}
^{n}W_{i_2},$ in $W=\sum_{j \in J}W_j,$ and continue in this
manner until no two eigenvalues for the $^s \nabla_{t_i
\partial/ \partial t_i}$ differ by a non-zero integer, and denote by
$\Lambda '$ the lattice generated by this. Then Th\'{e}or\`{e}me
3.4. in loc. cit. shows that $\Lambda '$ is free, invariant under
the $\nabla_{t_i\partial/ \partial t_i},$ and $\nabla$ is
represented on $\Lambda'$ as in the statement of the lemma. \hfill
$\Box$

\begin{lemma}
$\hat{\varphi}^*$ is an equivalence of categories.
\end{lemma}

{\it Proof.} In order to prove this, by the definition above, we
may assume that $\hat{\varphi}^*$ is  $i^*,$ where $i$ is the
closed immersion
$$
i:{\rm Spf} K[[t_1,\cdots, t_r ]]_{t_1 \cdots t_r} \to {\rm Spf}
K[[t_1,\cdots,t_n ]]_{t_1 \cdots t_r}.
$$
Since $i^*$ has  a section, namely $pr^{*}$, where $pr$ is the
projection, the essential surjectivity follows. To see the
full-faithfullness it is enough to show, for $$(E,\nabla) \in {\rm
Mic}_{reg,t_1 \cdots t_r }(K[[t_1,\cdots, t_n  ]]/K) ,$$ that
$$
\Gamma_{\nabla}(K[[t_1,\cdots, t_n  ]]_{t_1 \cdots t_r },E) \to
\Gamma_{i^* \nabla} (K[[t_1,\cdots, t_r ]]_{t_1 \cdots t_r } ,i^*
E)
$$
is an isomorphism. By lemma 3. we may assume that
$$(E,\nabla)=(\mathcal{O}^{\oplus rkE},d-\sum_{1 \leq i \leq r} \Gamma_i
\frac{dz_i}{z_i})$$ such that the only  integer eigenvalues of the
$\Gamma_{i}$ are 0. Then we see  that  $U(z):=\sum_{I \in
\mathbb{Z}^r}U_{I }z^I \in \Gamma_{\nabla}(K[[t_1,\cdots, t_n
]]_{t_1 \cdots t_r },E)$ if and only if $U_{0} \in \cap_{1\leq
i\leq r} \Gamma_{i}$ and $U_I=0,$ for $I \in \mathbb{Z}^r
\setminus \{0 \}.$ This shows the bijectivity. \hfill $\Box$

Finally in order to define the tangential basepoints we will need
the following
\begin{lemma}
Let ${\rm Mic}_{reg}(K[t_1,\cdots, t_n  ]_{t_1 \cdots t_n}/K)$
denote the category of vector bundles with integrable connection
and regular singularities at infinity on $(\mathbb{A}^{1}
\setminus \{ 0\})^n $. Then the natural map
$$
{\rm Mic}_{reg }(K[t_1,\cdots, t_n  ]_{t_1 \cdots t_n}/K) \to {\rm
Mic}_{reg,t_1 \cdots t_r }(K[[t_1,\cdots, t_n  ]]/K)
$$
is an equivalence of categories.
\end{lemma}

{\it Proof.} Since $(\mathcal{O}^m, d- \sum_{1\leq i\leq
n}\Gamma_{i}\frac{dz_i}{z_i}),$ with $\Gamma_{i} \in M_{m\times
m}(K) ,$ has regular singularities along all the divisors at
infinity, the essential surjectivity of the functor follows from
lemma 3.

Let $(E,\nabla) \in {\rm Mic}_{reg }(K[t_1,\cdots, t_n  ]_{t_1
\cdots t_n}/K)$ then since $E$ has a locally free extension to
$K[t_1,\cdots,t_n ],$ e.g. a logarithmic extension of
$(E,\nabla),$ the natural map
$$\Gamma_{\nabla}(K[t_1,\cdots,t_n]_{t_1 \cdots t_n },E) \to
\Gamma_{\nabla}(K[[t_1,\cdots,t_n]]_{t_1 \cdots t_n },E)
$$
is injective. To see its surjectivity, let $s$ be in the target.
Then  $(E,\nabla)$ and  $s$ can be defined over a subfield $K_0$
of $K$ that has countable transcendence degree over $\mathbb{Q}.$
Therefore it suffices to prove the surjectivity where $K$ is
replaced by  an arbitrary subfield $L$ of $\mathbb{C}.$ If
$L=\mathbb{C}$ then by the Riemann-Hilbert correspondence we may
assume that $(E,\nabla)= (\mathcal{O}^{\oplus k}, d- \sum_{1\leq i
\leq n}\Gamma_i \frac{dz_i}{z_i} ),$ where $\Gamma_i \in M_{k
\times k}(\mathbb{C})$ with 0 as the only possible integer
eigenvalue for the $\Gamma_i'$s. Then any formal solution of the
differential equation around 0 is  a constant section of
$\mathcal{O}^{k}$ that is annihilated by all the $\Gamma_i'$s.
This gives the surjectivity when $L=\mathbb{C}.$ For general $L
\subseteq \mathbb{C},$ let $s \in
\Gamma_{\nabla}(L[[t_1,\cdots,t_n]]_{t_1 \cdots t_n },E).$ Then,
by the above, there is an $s' \in
\Gamma_{\nabla}(\mathbb{C}[t_1,\cdots,t_n]_{t_1 \cdots t_n },E)$
such that $s'$ maps to the pull-back of $s$ by $\mathbb{C}/L.$
Then by the injectivity proven above, $s'$ is invariant under
${\rm Gal }(\mathbb{C}/L)$ since $s$ is. Therefore by descent $s$
comes from the variety over $L.$ \hfill $\Box$

Therefore we obtain a quasi-inverse
$$
{\rm Mic}_{reg,\hat{D}_{(x)}}(\hat{N}/K) \to {\rm Mic}_{reg
}(N_{D_x/\overline{X}} ^* /K),
$$
and combining this with
$$
{\rm Mic}_{reg}(X/K) \to {\rm Mic}_{reg,\hat{D}_{(x)}}(\hat{N}/K )
$$
and $\hat{\varphi} ^*$ gives
$$
{\rm Mic}_{reg}(X/K) \to {\rm Mic}_{reg }(N_{D_x/\overline{X}} ^*
/K)
$$

{\bf 5.2. Functoriality.} Let $f:\overline{X}_{log} \to
\overline{Y}_{log}$ be a morphism of log schemes with $f(x)=y.$
Let $\mathcal{P}(f): N_{x,X,log} \to N_{y,Y,log}$ denote the
unique homogeneous map that induces $f_{x}:x_{log}\to y_{log}$
under restriction and the isomorphisms $x_{log}\simeq 0_{x,log}$
and $y_{log}\simeq 0_{y,log}.$ $\mathcal{P}(f)$ is called the {\it
principal part } of $f$ at $x$ relative to the given divisors.
Then we have a commutative diagram
\[
\begin{CD}
{\rm Mic}_{reg,E}(Y/K) @>{f^*}>> {\rm Mic}_{reg,D}(X/K)\\
@VVV  @VVV\\
{\rm Mic}_{reg,E_{(y)}}(N_{y}/K) @>{\mathcal{P}(f)^*}>> {\rm
Mic}_{reg,D_{(x)}}(N_{x}/K).
\end{CD}
\]

\section{The unipotent case}

{\bf 6.1. } We will need the following in order to give the
relation between tangential and ordinary basepoints. In a
tannakian category we denote the tannakian subcategory of
unipotent objects by the subscript {\it uni}.
\begin{lemma}
Let $D_i \subseteq (\mathbb{P}_{K}^1)^r,$ for $i \in \{0,\infty
\}$ denote the divisor defined by the standard coordinate axes
passing through $i.$ And let $D:=D_0 \cup D_{\infty}.$ Then the
restriction to the origin map
$${\rm Mic}_{uni}((\mathbb{P}_{K}^1)^r_{log}/K) \to {\rm
Mic}_{uni}(K_{log}^r/K)
$$
is an equivalence of categories.
 \end{lemma}

{\it Proof.} Denote $(\mathbb{P}^{1})^r$ by $P.$ First note that
${\rm Mic}_{uni}(K^{r}_{log}/K)$ is canonically equivalent to
$\mathcal{T}_{uni,r}/K$ the full subcategory of unipotent objects
in $\mathcal{T}_r/K,$ where by $\mathcal{T}_r$ we mean the
category defined as in Definition. 1  above with ${\rm
Cart^{-}}(\cdot)$ replaced with $\mathbb{N}^r.$  Under this
equivalence the functor above becomes the one that associates
$E(0)$ endowed with $\{ {\rm res}_{i}(E,\nabla)\}_{1\leq i \leq r}
$ to $(E,\nabla),$ where ${\rm res}_{i}$ denotes the residue along
$D_i$ at $0.$

Let $(V, \{ N_{i}\}_{i \in I}) $ be an object of
$\mathcal{T}_{r,uni}$ then $(V \otimes_{K} \mathcal{O}_{P}, d-
\sum_{i \in I }N_i d \log z_i)$ is an object of ${\rm Mic}_{uni
}(P_K/K )$ that has image $(V, \{ N_{i}\}_{i \in I}) $ under the
restriction functor. This proves the essential surjectivity.

Since the functor above  is a tensor functor in order to see that
it is fully faithful it is enough to show that taking fibers at
zero induces an isomorphism
$$
{\rm Hom}_{\nabla}((\mathcal{O}_{P_K},d),(E,\nabla ) ) \to  {\rm
Hom }_{\mathcal{T}_r}((K,\{ 0 \}_{i \in I} ) , (E(0),{\rm \{
res_i(E,\nabla)} \}_{i \in I} )
$$
is an isomorphism, or equivalently that
$$
 {\rm H }^{0}_{dR}(P_{K,log}, (E,\nabla) ) \to \cap_{i \in I} {\rm
ker}_{E(0)}({\rm res}_{i}(E,\nabla))
$$
is an isomorphism. First note that for the underlying bundle $E$
of $$(E,\nabla) \in {\rm Mic}_{uni}(P_K/K)$$ is trivial. This
follows from the fact that ${\rm
Ext}^{1}_{P}(\mathcal{O}_P,\mathcal{O}_P)={\rm
H}^{1}(P,\mathcal{O}_{P})=0$ by induction on the nilpotence level
[De]. Therefore without loss of generality we will assume  that
$(E,\nabla)=(\mathcal{O}_P ^{rk E}, d- \sum_{1\leq i \leq
r}N_{i}\, d \log z_i), $  for some nilpotent matrices $N_{i} \in
M_{rkE \times rkE} (K). $ If  $\alpha $ is  a global (horizontal)
section of $(E,\nabla)$ then  it is a constant section of
$\mathcal{O}^{rk E}.$ Therefore the map above is injective. In
order to see that it is surjective,  we note that  for any $\alpha
\in \cap_{i \in I} {\rm ker}_{E(0)}({\rm res}_{i}(E,\nabla)),$ the
constant section of $\mathcal{O}_{P}^{rk E }$ with fiber $\alpha$
at 0 is a horizontal section with respect to the connection $d -
\sum_{i \in I}N_i d \log z_i.$ \hfill $\Box$

Let $\overline{N}_{D_x/\overline{X}}:=\prod_{i \in I
}\overline{N}_{D_{i}/\overline{X}},$ where
$\overline{N}_{D_i/\overline{X}}$ is the smooth compactification
of $N_{D_i/\overline{X}}.$  Let $\overline{D}_{(x)} \subseteq
\overline{N}_{D_x/\overline{X}}$ denote the normal crossings
divisor $D_0 \cup D_{\infty}.$

By the first section we have a functor
$$
{\rm Isoc}_{uni}(\overline{X}_{log}/W) \to {\rm
Isoc}_{uni}(x_{log}/W),
$$
combining this with the canonical isomorphism $x_{log}\simeq
0_{log},$ in section 3 and applying the last lemma to obtain an
equivalence of categories
$$
{\rm Isoc}_{uni}(\overline{N}_{log}/W) \to {\rm Isoc}_{uni}
(0_{log}/W ),$$ we obtain a functor
$$
\varphi^*:{\rm Isoc}_{uni}(\overline{X}_{log}/W) \to  {\rm
Isoc}_{uni}(\overline{N}_{log}/W ).
$$

{\rm Remark.} Note that in the de Rham version of the above we
have a commutative diagram
\[
\begin{CD}
{\rm Mic}_{uni}(\overline{X}_{log}/K) @>>>  {\rm
Mic}_{uni}(\overline{N}_{log}/K)\\
@VVV  @VVV \\
{\rm Mic}_{reg}(X/K) @>>> {\rm Mic}_{reg}(N/K).
\end{CD}
\]

{\bf 6.2. Description of $\varphi^*.$} First we will give another
description of the functor $${\rm
Isoc}_{uni}(\overline{X}_{log}/W) \to {\rm Isoc}_{uni}(0_{log}/W).
$$

   We use the notation of section
2.(iii)., i.e. $\overline{\mathfrak{X}}$ is a lifting of
$\overline{X}$ etc. Let $\psi:\mathcal{U}_{log} \to
\overline{\mathfrak{X}}_{log,K} ,$  where $\mathcal{U} \subseteq
(\mathfrak{N}_{\mathfrak{D}_{\mathfrak
x}/\overline{\mathfrak{X}}})_K$ is a polydisc around zero, be a
map such that $\psi(0)=\mathfrak{x}_K;$ the  map on the monoids is
the identity map; the map induced by the differential at zero,
$d_{0} (\psi)$ from $\oplus_{i \in I
}(\mathfrak{N}_{\mathfrak{D}_i/\overline{\mathfrak{X}}})_K \simeq
T_{0}(
\mathfrak{N}_{\mathfrak{D}_{\mathfrak{x}}/\overline{\mathfrak{X}}})_{K}$
to $\oplus_{i \in I
}(\mathfrak{N}_{\mathfrak{D}_i/\overline{\mathfrak{X}}})_K$ is the
identity map; and $\psi$ is an isomorphism onto its image. We have
a map
$$
\psi^*: {\rm Mic}_{uni}(\overline{\mathfrak{X}}_{ K, log}/K ) \to
{\rm Mic}_{uni}(\mathcal{U}_{log}/K),
$$
where  $\mathcal{U}_{log}$ is the   log scheme associated to the
divisor obtained by pulling back $\mathfrak{D}_{K}.$ Combining
this with the restriction map, we obtain
$$
\psi^* _0 : {\rm Mic}_{uni}(\overline{\mathfrak{X}}_{K,log}/K) \to
{\rm Mic}_{uni}(0_{\mathfrak{x},log}/K).
$$
Note that here $0_{\mathfrak{x},log}$ denotes  a log scheme over
$K$ and that in fact the category ${\rm
Mic}_{uni}(0_{\mathfrak{x},log}/K)$ does not depend, up to
canonical isomorphism, on the choices of the liftings since ${\rm
Mic}(0_{\mathfrak{x},log}/K)\simeq \mathcal{T}_{\mathfrak{x}}/K$
and  ${\rm Weil}^{-}(\mathcal{U},\mathcal{D}_{\mathfrak{x}})\simeq
{\rm Weil}^{-}(\overline{X},D_{x}).$

Let $\mathcal{V}\subseteq
(\mathfrak{N}_{\mathfrak{E}_{\mathfrak{y}}/\overline{\mathfrak{Y}}})_K,$
and $\tilde{\psi}:\mathcal{V} \to \overline{\mathfrak{Y}}_K $ be
another choice as above, with $\tilde{\psi}(0)=\mathfrak{y}_K.$
Then we have a map
$$
(\psi \times \tilde{\psi})^{\sim}: (\mathcal{U} \times
\mathcal{V})_{log}^{\sim} \to (\overline{\mathfrak{X}}_{K}\times
\overline{\mathfrak{Y}}_{K})_{log}^{\sim}
$$
induced by $\psi \times \tilde{\psi}.$ The underlying map of
schemes is the identity map on the exceptional divisors. Note that
the exceptional divisors are respectively the products of the
normal bundles of $\psi^*(\mathfrak{D}_{K,i}) \times \tilde{\psi}
^*(\mathfrak{D}_{K,i}) $ and $\mathfrak{D}_{K,i}\times
\mathfrak{D}_{K,i}$ at $(0,0)$ and
$(\mathfrak{x}_{K},\mathfrak{y}_{K}).$ If $(E,\nabla) \in {\rm
Isoc}_{uni}(\overline{X}_{log}/W)$ then by pulling back with
$(\psi\times \tilde{\psi})^{\sim}$ we have a canonical isomorphism
between the pullbacks of $\psi^*
_{0}(E,\nabla)_{\overline{\mathfrak{X}}}$ and $\tilde{\psi}^*
_0(E,\nabla)_{\overline{\mathfrak{Y}}}$ to the tube of the
diagonal in  $(0_{\mathfrak{x},log}\times
0_{\mathfrak{y},log})^{\sim}.$ Here note that even though the
diagonal is not defined, the tube of the diagonal is well-defined
as the tube of the diagonal after $0_{\mathfrak{x},log}/W$ and
$0_{\mathfrak{y},log}/W$ are identified by an isomorphism that
induces the identity map on the special fibers. These isomorphisms
on the tubes, viewed as log rigid analytic spaces, satisfy the
cocycle condition and hence define an element of ${\rm
Isoc}_{uni}(0_{log}/W)$  by lemma 5 and 6.

\begin{lemma}
Let $(\mathfrak{N},\mathfrak{D})/W$ and
$(\mathfrak{M},\mathfrak{E})/W$  be two formal vector bundles
endowed with linear normal crossings divisors passing through 0,
which are liftings of a vector bundle with normal crossings
divisor $(N,D)/k.$ Using the standard notation as above, let
$(E,\nabla)$ and $(F,\nabla)$ be in ${\rm
Mic}_{uni}(\overline{\mathfrak{N}}_{K,log}/K )$ and ${\rm
Mic}_{uni}(\overline{\mathfrak{M}}_{K,log}/K).$   And let
$$
\alpha_{0}: p_{1} ^*(E,\nabla)_0 \to p_2 ^* (F,\nabla)_0
$$
be  an isomorphism between the pullbacks to the tube of the
diagonal in  $(0_{log}\times 0_{log})^{\sim }$ of the restrictions
of $(E,\nabla)$ and $(F,\nabla)$ to $0_{log}$ in
$\mathfrak{N}_{K}$ and $\mathfrak{M}_{K}.$ Then there is a unique
isomorphism
$$
\alpha: p_{1} ^* (E,\nabla) \to p_2 ^* (F,\nabla)
$$
on the tube of the diagonal  in $(\overline{\mathfrak{N}}\times
\overline{\mathfrak{M}})^{\sim}_{K,log}$ that extends $\alpha_0.$
\end{lemma}

{\it Proof.} First we will assume without loss of generality that
the special fiber of the formal vector bundles is the trivial
bundle $k^r.$  In order to prove the lemma  we will assume,
without loss of generality by choosing bases whose reductions
modulo p are the standard basis of $k^{r},$ that
$(\mathfrak{N},\mathfrak{D})$ and $(\mathfrak{M},\mathfrak{E})$
are $\mathbb{A}^{r}$ endowed with the standard coordinate
hyperplanes as the divisor. Furthermore we will assume by Lemma 3.
that $(E,\nabla)=(F,\nabla) = (\mathcal{O}^{r}, d- \sum_{i \in I}
N_{i}\,d \log z_i),$ for some nilpotent operators $N_{i}.$ The
existence follows from the fact that unipotent logarithmic
connections on the generic fiber p-adically converge, that is
induce logarithmic isocrystals. More explicitly the isomorphism
from  $p_1 ^* (E,\nabla)$ to $p_2 ^* (E,\nabla)$ on the tube of
the diagonal is given by $ \prod_{1\leq i\leq r}
e^{\log(1+t_i)N_i} ,$ where we let $\frac{z_i ^{(2)}}{z_{i}
^{(1)}}=1+t_i,$ for $|t_i|<1.$ In particular note that
$\alpha((\mathfrak{v}_K
,\mathfrak{u}_{K}))=\alpha_0([\mathfrak{v}_K ,\mathfrak{u}_{K}]),$
when $E=V_{1} \otimes \mathcal{O}$ and $F=V_2 \otimes \mathcal{O}
$ as in the statement of the lemma, where $\mathfrak{u}$ and
$\mathfrak{v}$ have the same reduction $v \in N^{*}(k)$
 modulo p, and $[\mathfrak{v}_{K},\mathfrak{u}_K]$ denotes the
 point corresponding to $(\mathfrak{v}_K,\mathfrak{u}_{K})$
 in the exceptional divisor.

In order to see the uniqueness it is enough to show that for a
unipotent vector bundle with logarithmic connection $(E,\nabla)$
on the tube
$]\Delta_{\overline{N}}[_{(\overline{\mathfrak{N}}\times
\overline{\mathfrak{M}})^{\sim}}$ of the diagonal in
$(\overline{\mathfrak{N}} \times \overline{\mathfrak{M}})^{\sim} $
the restriction functor
$$
{\rm
H}^{0}_{dR}(]\Delta_{\overline{N}}[_{(\overline{\mathfrak{N}}\times
\overline{\mathfrak{M}})^{\sim}},(E,\nabla)) \to {\rm
H}^{0}_{dR}(]0[_{(0_{log}\times 0_{log})^{\sim}},(E,\nabla))
$$
is injective. This follows from ${\rm
H}^{0}_{dR}(]\Delta_{\overline{N}}[_{(\overline{\mathfrak{N}}\times
\overline{\mathfrak{M}})^{\sim}},(\mathcal{O},d))=K$ by induction
on the nilpotence level of $(E,\nabla).$ \hfill $\Box$

Therefore we obtain a map
 $${\rm Isoc}_{uni}(\overline{X}_{log}/W) \to
{\rm Isoc}_{uni}(0_{\mathfrak{x},log}/W). $$ In order to see that
this is the same as the one constructed before we note that  the
map
$$
(0_{\mathfrak{x},log}\times 0_{\mathfrak{y},log})^{\sim} \to
(\mathcal{U} \times \mathcal{V})_{log}^{\sim} \to
(\overline{\mathfrak{X}}_{K}\times
\overline{\mathfrak{Y}}_{K})_{log}^{\sim}
$$
induced by $\psi$ and $\tilde{\psi}$ factors as
$$
(0_{\mathfrak{x},log}\times 0_{\mathfrak{y},log})^{\sim} \to
(\mathfrak{x}_{K,log} \times \mathfrak{y}_{K,log})^{\sim} \to
(\overline{\mathfrak{X}}_{K}\times
\overline{\mathfrak{Y}}_{K})_{log}^{\sim}.
$$
The underlying schemes of  $(0_{log} \times 0_{log})^{\sim }$ and
$(\mathfrak{x}_{K,log} \times \mathfrak{y}_{K,log})^{\sim}$ are
the exceptional divisors of $(\mathcal{U} \times
\mathcal{V})_{log}^{\sim},$ and
$(\overline{\mathfrak{X}}_{K}\times
\overline{\mathfrak{Y}}_{K})_{log}^{\sim}$ and the map induced
between them is the identity map if they are identified via the
isomorphisms of the form $\mathfrak{z}_{K,log}\simeq
0_{\mathfrak{z},log},$ $\mathfrak{z} \in
\{\mathfrak{x},\mathfrak{y} \}$ because of the condition on the
derivative at 0 of the local isomorphisms $\psi$ and
$\tilde{\psi}.$ This shows that the two functors from ${\rm
Isoc}_{uni}(\overline{X}_{log}/W)$ to ${\rm
Isoc}_{uni}(0_{\mathfrak{x},log}/W)$ are the same.

To give a description of $\varphi^*$ we need to describe
$$
{\rm Isoc}_{uni}(0_{log}/W) \to {\rm
Isoc}_{uni}(\overline{N}_{log}/W).
$$
Let
$\mathfrak{N}_{\mathfrak{D}_{\mathfrak{x}}/\overline{\mathfrak{X}}}/W$
be a formal vector bundle, lifting $N_{D_x/\overline{X}}$ endowed
with a linear normal crossings divisor at 0. We do not assume that
the lifting in fact comes from the tangent space of a lifting
$\overline{\mathfrak{X}}$ of $\overline{X}$ but continue to  use
that  notation for consistency. The functor,  defined only up to
canonical isomorphism,
$$
e_{\overline{\mathfrak{X}}}: \mathcal{T}_{x,uni}/K \simeq {\rm
Mic}_{uni}(0_{\mathfrak{x},log}/K) \to {\rm
Mic}_{uni}((\overline{\mathfrak{N}}_{\mathfrak{D}_{\mathfrak{x}} /
\overline{ \mathfrak{X}}})_{K,log}/K)
$$
is the one that sends $(V, \{N_i \}_{i \in I})$ to $ (V \otimes _K
\mathcal{O}_{\overline{\mathfrak{N}}},d -\sum_{i \in I} N_i d \log
z_i ),$ where $\{ z_i\}_{i \in I}$ is a linear system of
coordinates for $\mathfrak{N}_K$ such that the divisor on
$\mathfrak{N}_{K}$ is defined by $\prod_{i \in I}z_i=0.$ Note that
since the $z_i$ are defined only up to a scalar multiple the $d
\log z_i$ are well-defined.

Let
$\mathfrak{N}_{\mathfrak{E}_{\mathfrak{y}}/\overline{\mathfrak{Y}}}/W$
be another such lifting. Then for $(E,\nabla) \in {\rm
Isoc}_{uni}(0_{log}/W),$ by the construction we have
$$e_{\overline{\mathfrak{X}}}(E,\nabla) \in {\rm
Mic}_{uni}((\overline{\mathfrak{N}}_{\mathfrak{D}_{\mathfrak{x}}/\overline{\mathfrak{X}}
})_{K,log} /K) {\rm \;\;and \;\;} e_{\overline{\mathfrak{Y}}}
(E,\nabla) \in
 {\rm
Mic}_{uni}((\overline{\mathfrak{N}}_{\mathfrak{E}_{\mathfrak{y}}/\overline{\mathfrak{Y}}
})_{K,log} /K)$$ and an isomorphism
$$
\alpha_{0}: p_{1} ^* e_{\overline{\mathfrak{X}}}(E,\nabla)_{0} \to
p_2 ^* e_{\overline{\mathfrak{Y}}}(E,\nabla)_{0}.
$$
And hence an isomorphism
$$
\alpha: p_1 ^* e_{\overline{\mathfrak{X}}}(E,\nabla) \to p_2 ^*
e_{\overline{\mathfrak{Y}}}(E,\nabla)
$$
by the last lemma. Since the isomorphisms on the log points
satisfy the cocycle condition, the isomorphisms for the
connections on the $\overline{\mathfrak{N}}_{K}'$s also satisfy
the cocycle conditon by the uniqueness statement in the last
lemma. Therefore the data of objects
$e_{\overline{\mathfrak{X}}}(E,\nabla),$ for each lifting
$\overline{\mathfrak{X}}$  as above, together with the
isomorphisms of their pullbacks to the blow-up of the products for
different liftings  define an element $e(E,\nabla) \in {\rm
Isoc}_{uni}(\overline{N}_{log}/W).$ This is a more explicit
description of the  inverse of
$$
{\rm Isoc}_{uni}(\overline{N}_{log}/W) \to {\rm
Isoc}_{uni}(0_{log}/W)
 $$
that will be useful below.

{\bf 6.3. Comparison with ordinary basepoints.} Let $Y/k$ be a
smooth variety, $y \in Y(k),$ $i: y \to Y$ the inclusion. Then we
denote the fiber functor $i^{*}: {\rm Isoc}(Y/W) \to Vec_K$ by
$\omega(y).$

Choosing $v \in N_{D_{x}/\overline{X}}^*(k),$ and composing
$\varphi^*$ with $$\omega(v):{\rm
Isoc}_{uni}(\overline{N}_{log}/W) \to {\rm Isoc}_{uni}(N^*/W) \to
Vec_{K}$$ we obtain the fiber functor $ \varphi_v ^*: {\rm
Isoc}_{uni}(\overline{X}_{log}/W) \to Vec_K .$

\begin{proposition} There is a canonical natural isomorphism $\varphi^*_v
\simeq fib(v).$ Here $fib(v)$ is the tangential basepoint defined
above.
\end {proposition}

{\it Proof.} Let $(E,\nabla) \in {\rm
Isoc}_{uni}(\overline{X}_{log}/W)$ and let
$\overline{\mathfrak{X}},$ etc. be a lifting, and let $\psi:
\mathcal{U} \to \overline{\mathfrak{X}}_{K}$ be  a map as above.
Then by the description of $\varphi^*$ above, $$(\varphi^*
(E,\nabla))_{\overline{\mathfrak{X}},\varphi}=
(E_{\overline{\mathfrak{X}}}(\mathfrak{x}_{K})\otimes_{K}
\mathcal{O}_{\overline{\mathfrak{N}} _{K}} ,  d - \sum_{i \in I}
{\rm
res}_{\mathfrak{D}_{K,i}}(E,\nabla)_{\overline{\mathfrak{X}}}(\mathfrak{x}_{K})
d \log z_i)
$$
since ${\rm res}_{\varphi^{-1}( \mathfrak{D}_{K,i} )}
\varphi^{-1}(E,\nabla)_{\overline{\mathfrak{X}}} (\mathfrak{x}_K)=
{\rm res}_{ \mathfrak{D}_{K,i} }
(E,\nabla)_{\overline{\mathfrak{X}}}(0).$ Together with the
choices above let $\mathfrak{v} \in
\mathfrak{N}_{\mathfrak{D}_{\mathfrak{x}} /\overline{\mathfrak{X}}
}^*(W)$ lifting $v.$ Then by the above formula, the realization
$\varphi^{*}_{v}(\overline{\mathfrak{X}}, \psi, \mathfrak{v} )$ of
$\varphi_{v} ^*$ is the map sending $(E,\nabla)$ as above to
$(E_{\overline{\mathfrak{X}}}(\mathfrak{x}_K)\otimes _{K}
\mathcal{O}_{\overline{\mathfrak{N}}_{K}
})(\mathfrak{v}_K)=E_{\overline{\mathfrak{X}}}(\mathfrak{x}_K) .$
Similarly by the description given in 1.(iii).
$fib(v)_{\overline{\mathfrak{X}},\mathfrak{v}}
(E,\nabla)=E_{\overline{\mathfrak{X}}}(\mathfrak{x}_K).$

Let $(\overline{\mathfrak{Y}}, \tilde{\psi}, \mathfrak{u}) $ be
another such data of a lifting. Then the isomorphism between
$fib_{\overline{\mathfrak{X}},\mathfrak{v}} (E,\nabla)$ and
$fib_{\overline{\mathfrak{Y}},\mathfrak{u}} (E,\nabla)$ is the one
obtained from pulling the connections to the point
$[\mathfrak{v}_{K},\mathfrak{u}_{K}] \in
(\overline{\mathfrak{X}}_{K} \times
\overline{\mathfrak{Y}}_{K})^{\sim}.$  Similarly the isomorphism
between $\varphi_{v} ^* (\overline{\mathfrak{X}},\psi,\mathfrak{v}
) (E,\nabla)$ and $\varphi_{v} ^*
(\overline{\mathfrak{Y}},\tilde{\psi},\mathfrak{u} ) (E,\nabla)$
is the one obtained by pulling back the connections to the point
$(\mathfrak{v}_{K},\mathfrak{u}_{K}) \in
(\overline{\mathfrak{N}}_{\overline{\mathfrak{X}}}\times
\overline{\mathfrak{N}}_{\overline{\mathfrak{Y}}})^{\sim} _{K}.$
Therefore in order to prove the claim we only need to see that for
$(V_{1} \otimes \mathcal{O} ,\nabla) \in {\rm
Mic}_{uni}((\overline{\mathfrak{N}}_{\overline{\mathfrak{X}}})_{K,log}/K)$
 and $(V_{2} \otimes \mathcal{O} ,\nabla) \in {\rm
Mic}_{uni}((\overline{\mathfrak{N}}_{\overline{\mathfrak{Y}}})_{K,log}/K)$
and an isomorphism between their pullbacks to
$(\overline{\mathfrak{N}}_{\overline{\mathfrak{X}}}\times
\overline{\mathfrak{N}}_{\overline{\mathfrak{Y}}})^{\sim}_{K},$
the isomorphism at the point $[\mathfrak{v},\mathfrak{u}]$ is the
same as the one at $(\mathfrak{v},\mathfrak{u} ).$ But this was
shown in the proof of Lemma 4. \hfill $\Box$

{\bf 6.4. Tangential basepoints for unipotent overconvergent
isocrystals.} We continue to use the notation above, and let ${\rm
Isoc}_{\overline{X}}^{\dagger}(X/W)$ denote the category of
isocrystals overconvergent along $D\subseteq \overline{X}.$ Then
the natural map
$$
{\rm Isoc}_{uni}(\overline{X}_{log}/W) \to {\rm
Isoc}_{\overline{X},uni}^{\dagger}(X/W)
 $$
is an equivalence of categories. This follows, in the usual way,
from the fact that for any log $F^{k}$-isocystal $\mathcal{E},$ in
particular $\mathcal{O}_{\overline{X}}\otimes K,$ on
$\overline{X}_{log}$ the natural map
$$
{\rm H}_{rig}^{i}(\overline{X}_{log}/W,\mathcal{E}) \to {\rm
H}^{i}_{rig,\overline{X}}(X/W, j^{\dagger}\mathcal{E})
$$
 is an isomorphism  ([Shi], proof of (2.4.1)) for all $i,$ where $j: X
\hookrightarrow \overline{X} $ is the inclusion. Therefore by the
above construction we obtain a fiber functor on ${\rm
Isoc}_{\overline{X},uni}^{\dagger}(X/W),$ and hence a fiber
functor
$$
fib(v): {\rm Isoc}_{uni}^{\dagger}(X/W) \to Vec_{K},
$$ for $v \in  N^{*}_{D_x/\overline{X}}(k).$

\section{Drinfel'd-Ihara relation}

{\bf 7.1. de Rham basepoint.} From now on let $X:=\mathbb{P}^1
\setminus \{ 0, 1, \infty\},$ and $\overline{X}:=\mathbb{P}^1.$ In
the following when we write $\mathbb{P}^1/\mathbb{Q}$ we will
always assume that it is endowed with a specific choice of a
coordinate function, i.e.  a rational function $z\in
k(\mathbb{P}^1 /\mathbb{Q})$ such that
$\mathbb{Q}(z)=k(\mathbb{P}^1/\mathbb{Q}),$ where
$k(\mathbb{P}^1/\mathbb{Q} )$ denotes the field of rational
functions on $\mathbb{P}^1/\mathbb{Q}.$

Let
$$M_{0,n}:=\{(x_0, \cdots, x_{n-1})\in (\mathbb{P}^1) ^n| { \rm \;  all \; } x_i {\rm \; are \; different} \}/PGL_2,$$
where $PGL_2$ acts diagonally by linear fractional
transformations.

\begin{lemma}
For $n \geq 4,$ $M_{0,n}$ has a compactification
$\overline{M}_{0,n}$ with  $\overline{M}_{0,n} \setminus M_{0,n}$
a simple normal crossings divisor and
$H^{1}_{B}((\overline{M}_{0,n})_{\mathbb{C}},\mathbb{Z})=0.$
\end{lemma}

{\it Proof.} First note that by using a linear fractional
transfornation that sends $(x_0,x_1,x_2)$ to $(0,1,\infty)$ we can
identify $M_{0,n}$ with $X^{n-3} \setminus \Delta_{h},$ where
$$\Delta_h:= \{(x_1,\cdots,x_{n-3}) | x_i \neq x_j {\rm \; for \;}
i \neq j \}$$ is the hyperdiagonal in $X^{n-3}.$   A
compactification $\overline{M}_{0,n}$ of $M_{0,n} \subseteq
(\mathbb{P}^{1})^{n-3}$ with $\overline{M}_{0,n} \setminus
M_{0,n}$ a simple normal crossings divisor is obtained by a
succession of blowings up of $(\mathbb{P}^{1})^{n-3}$ along linear
subvarieties. Since $(\mathbb{P}^{1})^{n-3}$ is simply connected
and the blowings up do not change the first cohomology group we
have $H^{1}_{B}((\overline{M}_{0,n})_{\mathbb{C}},\mathbb{Z})=0. $
\hfill $\Box$

The above lemma together with Grothedieck's comparison theorem
gives $$H^{1}
_{dR}((\overline{M}_{0,n})_{\mathbb{C}},(\mathcal{O},d))=0,$$ and
hence $H^{1}(\overline{M}_{0,n},\mathcal{O})=0.$ Therefore one has
a canonical fiber functor on $${\rm
Mic}_{uni}(M_{0,n}/\mathbb{Q})$$ defined in [De].  This can be
described as follows. For any  $(E,\nabla)\in {\rm
Mic}_{uni}(M_{0,n}/\mathbb{ Q}),$ the underlying vector bundle
$\overline{E}$ of the canonical extension $(\overline{E},\nabla)
\in {\rm Mic}_{uni}(\overline{M}_{0,n}, {\rm log} D )$ is trivial
[De], where $D:=\overline{M}_{0,n}\setminus M_{0,n}.$  Therefore
the canonical map $$\Gamma(\overline{M}_{0,n},\overline{E})
\otimes\mathcal{O}_{\overline{M}_{0,n}} \simeq \overline{E}$$ is
an isomorphism and the functor
$$
\omega(dR) \,(:=\Gamma(\overline{M}_{0,n},\cdot )): {\rm
Mic}_{uni}(M_{0,n}/\mathbb{Q}) \to {\rm Vec}_{\mathbb{Q}}
$$
 is a fiber functor.

This a priori depends on the choice of a compactification of
$M_{0,n}$ with zero first Betti cohomology. Let $\overline{M}^{1}
_{0,n}$ and $\overline{M}^{2} _{0,n}$ be two such
compactifications. By applying Hironaka's resolution of
singularities to the closure $\Delta(M_{0,n})^{-}$ of the image of
$M_{0,n}$ in $\overline{M}^{1} _{0,n}\times \overline{M}^{2} _{
0,n}$ under the diagonal map we obtain a compactification
$\overline{M}^{12} _{0,n}$ of $M_{0,n}$ with the complement a
simple normal crossings divisor together with maps
$\overline{M}^{12} _{0,n} \to \overline{M}^{1} _{0,n}$ and
 $\overline{M}^{12} _{0,n} \to \overline{M}^{2} _{0,n}$ that commute with the inclusions of $M_{0,n}.$  Let $(E,\nabla) \in {\rm Mic}_{uni}(M_{0,n} /\mathbb{Q}),$ and $(\overline{E}^{i},\nabla) $ its canonical extension to $\overline{M}^{i} _{0,n},$ for $i \in \{ 1,2\}.$   Then since ${\rm supp}(\pi_{i} ^* (D_i)) \subseteq D_{12}$ we have $\pi_i ^{*} (\overline{E} ^{i}\nabla ) \in {\rm Mic}(\overline{M}^{12} _{0,n} ({\rm log}D_{12}) /\mathbb{Q}).$ Furthermore since the exponents of the pull-backs are linear combinations of the original exponents, which are zero by the definition of canonical extension, the exponents of the pull-backs are zero as well. And hence $\pi_i ^{*} (\overline{E} ^{i}\nabla ) \in {\rm Mic}_{uni}(\overline{M}^{12} _{0,n} ({\rm log}D_{12}) /\mathbb{Q}).$ Therefore the  identification $\pi_1 ^{*} (\overline{E} ^{1}\nabla )|_{M_{0,n}} =\pi_2 ^{*} (\overline{E} ^{2}\nabla )$ extends to $\overline{M}^{12} _{0,n}$ and after taking global sections induces an isomorphism $\Gamma(\overline{M}^{12} _{0,n},\overline{E} ^{1}) \simeq \Gamma(\overline{M}^{12} _{0,n},\overline{E} ^{2}).$ A similar argument shows that the cocycle condition is satisfied for three different compactifications and hence the definition of the canonical de Rham fiber functor is independent of the compactification satisfying the properties above, and is functorial with respect to arbitrary maps between $M_{0,n}.$

If $x$ is a (tangential) basepoint of $M_{0,n}$ then the natural
maps $\Gamma(\overline{M}_{0,n},\overline{E}) \to \overline{E}(x)$
induce an isomorphism of the fiber functors $\omega(dR)$ and
$\omega(x).$ And hence there is a canonical path between the
fiber functors $\omega(x)$ and $\omega(y)$ in the de Rham theory;
this is denoted by $_{y} e(dR)_x.$ As above this path is
independent of the compactification in the case where $x$  and $y$
are ordinary points.

{\bf 7.2. Fundamental group of $M_{0,n}.$} Let
$\pi_{1,dR}(M_{0,n},\omega(dR))$ denote the fundamental group of
the tannakian category ${\rm Mic}_{uni}(M_{0,n}/\mathbb{Q})$ at
the fiber functor $\omega(dR);$ this is a pro-unipotent algebraic
group.

\begin{lemma}
The natural map $p_n: M_{0,n} \to M_{0,n-1}$ with  $$p_n([x_0,
\cdots,x_{n-1}])=[x_0, \cdots, x_{n-2}]$$ induces an exact
sequence
$$
0 \to \pi_{1,dR}(F_x, \omega(dR) ) \to
\pi_{1,dR}(M_{0,n},\omega(dR)) \to
\pi_{1,dR}(M_{0,n-1},\omega(dR))\to 0
$$
where $F_{x}$ is the fiber of $p_n$ containing $x:= [x_0,
\cdots,x_{n-1}].$
\end{lemma}

Proof. First assume that  the basefield is $\mathbb{C}.$ Then
$p_{n}$ is a locally trivial fibration with fibers isomorphic to
$\mathbb{P}^1$ minus $n-1$ points. Hence we get a homotopy exact
sequence for the topological fundamental groups
$$
\cdots \to \pi_{2}(M_{0,n-1},p_{n}(x) ) \to \pi_{1}(F_x,x) \to
\pi_{1}(M_{0,n},x ) \to \pi_{1}(M_{0,n-1},p_{n}(x)) \to \cdots.
$$
Since $F_x$ is connected  we obtain the exact sequence
$$
 \pi_{1}(F_x,x) \to \pi_{1}(M_{0,n},x) \to
\pi_{1}(M_{0,n-1},p_n(x) )\to 0
$$

Let $G$ be an abstract group, and $\{Z^{i}(G) \}_{1 \leq i}$ the
central descending series, i.e. $Z^{1}(G):=G$ and
$Z^{i+1}(G):=[G,Z^{i}(G)] $. Then $G^{[N]}:=(G/Z^{N+1}(G))/tors$
(note that since $G/Z^{N+1}(G)$ is nilpotent the set of its
torsion  elements form a  subgroup [Ba]) is canonically imbedded
into $G^{[N]}\, ^{*},$ the universal nilpotent torsion free
divisible group receiving a map from $G^{[N]}.$ If $g \in G^{[N]}
\, ^*$ then there is an $n \in \mathbb{Z}_{>0}$ such that $g^{n}
\in  G^{[N]}, $ and for every $g \in G^{[N]}$ and $n \in
\mathbb{Z}_{ \geq 0}$ there exists a unique $h \in G^{[N]} \,
^{*}$ such that $h^n=g$  ([Ba], section 8.3). These imply that
$$
\pi_{1}(F_{x},x)^{[N]}\, ^{*} \to \pi_{1}(M_{0,n},x)^{[N]}\, ^*
\to \pi_{1}(M_{0,n-1},p_n(x)  )^{[N]}\, ^* \to 0
$$
is exact.

For a pro-unipotent algebraic group $H,$ let $H^{(N)}$ denote its
largest quotient of nilpotence level $\leq N.$     Since, by the
Riemann-Hilbert correspondence, $\pi_{1,dR}(\cdot)^{(N)}$ is the
unipotent algebraic envelope of $\pi_{1}(\cdot)^{[N]*}$ over
$\mathbb{C},$  the exact sequence above implies  that
$$
\pi_{1,dR}(F_{x},x)^{(N)} \to \pi_{1,dR}(M_{0,n},x)^{(N)} \to
\pi_{1,dR}(M_{0,n-1},p_n(x) )^{(N)} \to 0
$$
is exact. Since any unipotent integrable connection on $F_x$ of
level $N$ can be extended to a unipotent integrable connection of
level $N$ on $M_{0,n}$ the first map is injective. By taking the
inverse limits we obtain the exact sequence
$$
0 \to \pi_{1,dR}(F_x,x) \to \pi_{1,dR}(M_{0,n},x ) \to
\pi_{1,dR}(M_{0,n-1},p_{n}(x) ) \to 0
$$
since the  inverse system $\{ \pi_{1,dR}(.\, ,\, .)^{(N)} \}_{N
\geq 1}$ satisfies the Mittag-Leffler condition. This gives the
statement in the lemma when the basefield is $\mathbb{C}.$ To
obtain it in the case when the base field is $\mathbb{Q}$ we note
that the unipotent de Rham fundamental group commutes with the
base change of the fields [De]. \hfill $\Box$

{\it Remark.} Note that the sequence
$$
0 \to \pi_{1}(F_x,x) \to \pi_{1}(M_{0,n},x) \to \pi_{1}
(M_{0,n-1},p_{n}(x)) \to 0
$$
is exact on the left as well. This follows from the fact that
$\pi_{i}(M_{0,n}, \; .)=0$ for $i \geq 2$ and $n \geq 4. $ To see
this  first note that the fibers $F^{(n)}$ of the maps $p_n,$
being isomorphic to $\mathbb{P}^{1}$ minus n-1 points have the
unit disc as a covering space (by uniformization theory),
therefore $\pi_{i}(F^{(n)}, \; . )=0$ for $i \geq 2.$ Then the
claim follows by induction from this and the homotopy sequences
for the fibrations $p_n.$

{\it Residues.} Let
$$
M^* _{0,n}  :=\{(x_0,\cdots, x_{n-1}) \in (\mathbb{P}^1)^{n} |
{\rm \; at \; most\; two\; of \; the \; } x_i {\rm \; are \;
equal}  \}/PGL_2.
$$
We have $M_{0,n} \subseteq M^* _{0,n}$ as an open subvariety with
the complement a simple normal crossings divisor. Let $D_{ij}
\subseteq M^* _{0,n}$ denote the divisor defined by the image of
$x_i=x_j.$  By the above lemma, using induction, we see that
$H_{1,B}(M_{0,n})$ is generated by the loops around the divisors
$D_{ij}.$ Therefore the image of $H_{1,B}(M _{0,n})$ in
$H_{1,B}(M^* _{0,n})$ is zero. Since $H_{1,B}(M^*
_{0,n},M_{0,n})=0,$ by using Mayer-Vietoris sequence and excision
for a suitable cover,  we see that $H_{1,B}(M^* _{0,n})=0.$
Therefore we may choose $\overline{M}_{0,n}$ as in the previous
lemma such that $M^* _{0,n}=\overline{M}_{0,n} \setminus D_{0},$
where $D_{0} \subseteq \overline{M}_{0,n}$ is a divisor contained
in $\overline{M}_{0,n}\setminus M_{0,n}.$

Let  $\overline{D}_{ij}$ denote the  closure of $D_{ij}$ in
$\overline{M}_{0,n}.$ Let $x \in
\overline{D}_{ij}(\overline{\mathbb{Q}})$ and $\{t_{k} \}_{1\leq
k\leq n-3} $ a system of parameters on $\overline{M} _{0,n}$ at
$x$ such that the divisor $\overline{M} _{0,n} \setminus M_{0,n}$
is defined by $t_1 \cdots t_r=0$ at $x$ and $\overline{D}_{ij}$
is defined by $t_1=0$ at $x.$ If we  let $(\overline{E},\nabla)$
denote the canonical extension of $(E,\nabla)$ we obtain a map
$$
-\nabla_{t_1 \frac{\partial}{\partial t_1}}(x): \overline{E}(x)\to
\overline{E}(x),
$$
which is independent of the choice of the system of parameters as
above. As usual we denote this map by ${\rm
res}_{x,\overline{D}_{ij}}(\overline{E},\nabla)$ This map
satisfies $$ {\rm res}_{x,\overline{D}_{ij}}((\overline{E},\nabla)
\otimes (\overline{F},\nabla) )= id_{\overline{E}(x)}\otimes  {\rm
res}_{x,\overline{D}_{ij}}((\overline{F},\nabla) + {\rm
res}_{x,\overline{D}_{ij}}((\overline{E},\nabla) \otimes
id_{\overline{F}(x)},
$$
and hence gives an element ${\rm res}_{x, \overline{D}_{ij}} \in
{\rm Lie} \pi _{1,dR}(M_{0,n}, x).$ Using the canonical
isomorphism
$$\Gamma((\overline{M}_{0,n})_{\overline{\mathbb{Q}}},\overline{E}
) \simeq \overline{E}(x)$$   we obtain an algebraic map
$$
(\overline{D}_{ij})_{\overline{\mathbb{Q}}} \, \to {\rm
End}(\Gamma(\overline{M}_{0,n})_{\overline{\mathbb{Q}}},
\overline{E}).
$$
Since $(\overline{D}_{ij})_{\overline{\mathbb{Q}}}$ is proper and
integral this map is in fact constant. And by descent this is in
fact defined over $\mathbb{Q}$  and we obtain  elements ${\rm
res}_{ij} \in {\rm Lie}\pi_{1,dR}(M_{0,n},\omega(dR)).$

If $\{ i,j\} \cap \{k,l \}=\emptyset $ then $D_{ij}\cap D_{kl}\neq
\emptyset $ and computing the residues at a point in $D_{ij}\cap
D_{kl}$ we see that $[{\rm  res}_{ij},{\rm res}_{kl}]=0$ by the
integrability of the connection. Similarly  by computing the
residues along a fiber of the projection $M_{0,n} \to M_{0,n-1}$
that maps $[x_0,\cdots,x_j, \cdots,x_{n-1}] $ to
$[x_0,\cdots,\hat{x}_j,\cdots,x_{x_n-1}],$ we see that
$\sum_{i}{\rm res}_{ij}=0,$ where we let ${\rm res}_{ii}=0$ by
convention. Let
$$
H_{n}:={\rm Lie}<< e_{ij}>>_{0\leq i,j \leq
n-1}/(e_{ii},e_{ij}-e_{ji},\sum_{i}e_{ij}, [e_{ij},e_{kl}]|{\rm \;
for \;} \{ i,j \} \cap \{ k,l\}=\emptyset ),
$$
with ${\rm Lie}<< \cdot>>$ denoting the free pro-nilpotent Lie
algebra generated by the arguments. We have a map $H_{n}\to {\rm
Lie }\pi_{1,dR}(M_{0,n},\omega(dR) ).$ This map is surjective
since for a space $Y$ that has a  compactification $\overline{Y}$
as above with
 $H^{1}_{B}(\overline{Y}_{\mathbb{C}})=0,$ ${\rm Lie
}\pi_{1,dR}(Y,\omega(dR))$ is generated by the dual of
$H^{1}_{dR}(Y/\mathbb{Q})$ [De].

\begin{corollary} The exact sequence in the statement of the previous lemma
has a natural splitting and gives
$$
\pi_{1,dR}(M_{0,n},\omega(dR)) \simeq \pi_{1,dR}(M_{0,n-1},
\omega(dR)) \ltimes \pi_{1,dR}(F^{(n)},\omega(dR)),
$$
where $F^{(n)} \simeq \mathbb{P}^1\setminus \{a_1,\cdots,a_{n-1}
\}$ for some dictinct points $a_i.$
\end{corollary}

Proof. Let
$$
H_{n}:={\rm Lie} <<e_{ij}>>_{1 \leq i,j \leq  n}/(e_{ii}, \,
e_{ij}-e_{ji}, \, \sum_{j}e_{ij}, \, [e_{ij},e_{kl}] |{\rm \, for
\, } \{i,j \} \cap \{ k, l\}=\emptyset   ).
$$

There are exact sequences
$$
0 \to G_{n} \to H_{n} \to H_{n-1} \to 0
$$
where $G_{n}:=Lie<<e_{in}>>/( \sum_{i}e_{in})  $

And there is a natural splitting $H_{n-1} \to H_n$ mapping
$e_{ij}$ to $e_{ij}.$

In other words
$$
H_n \simeq H_{n-1} \ltimes G_n
$$
where the action of $H_{n-1}$ on $G_n$ is determined by
\begin{eqnarray*}
[e_{ij},e_{kn}]&=&0 {\rm \; if \; } \{ i, j \} \cap \{ k, n
\}=\emptyset \\  {\rm and \; \;  }  [e_{ij}, e_{in}]&=&-[\sum_{k
\neq i} e_{kj}, e_{in}]=-[e_{jn},e_{in}] { \rm \; if \; } 1 \leq
i,j \leq n-1.
\end{eqnarray*}

As we have seen above there are  surjections
$$
H_{n} \to  {\rm Lie} \pi_{1,dR}(M_{0,n}, \omega(dR)).
$$
These fit into  commutative diagrams

\[
\begin{CD}
0 & \to  & G_n    & \to  & H_n & \to & H_{n-1} & \to & 0 \\
  & &  @VVV  @VVV    @VVV   \\
0 & \to & {\rm Lie}\pi_{1,dR}(F^{(n)}, \omega(dR)) & \to & {\rm
Lie}\pi_{1,dR}(M_{0,n},\omega(dR)) & \to & {\rm
Lie}\pi_{1,dR}(M_{0,n-1},\omega(dR)) &  \to & 0
\end{CD}
\]
with exact rows and this implies that they are in fact
isomorphisms
$$
H_{n} \simeq {\rm Lie} \pi_{1,dR}(M_{0,n}, \omega(dR))
$$
by induction.  \hfill $\Box$

From now on we will write $e_{ij}$ instead of ${\rm res}_{ij}.$

{\bf 7.3. Basepoints on $M_{0,4}$ and $M_{0,5}.$} For $i,j \in \{
0,1,\infty \}$ let $t_{ij}$ denote the unit tangent vector at the
point $i$  that points in the direction from $i$ to $j.$ For
example, $t_{01}:=\frac{d}{dz}$ at 0, $t_{10}:=-\frac{d}{dz}$ at
1, $t_{\infty 0}:=z^2 \frac{d}{dz}$ at $\infty$ etc. Note that
$X\simeq M_{0,4}$ by the map $z \to [0,z,1,\infty],$  and
$\mathbb{P}^{1}$ can be viewed as the compactification
$$
\{(x_0,x_1,x_2,x_3) | x_0 \neq x_2, \; x_0 \neq x_3, \; x_2 \neq
x_3 \}/PGL_2
$$
of $M_{0,4}.$ Therefore we may view  $t_{ij}$ as basepoints on
$M_{0,4}.$

Note that $M^* _{0,5}$ is a compactification of $M_{0,5}.$ We will
use the following tangential basepoints on $M^* _{0,5}$ at  the
points
\begin{eqnarray*}
& &[\{x_0,x_1\},\{x_2,x_3\},x_4],\; [\{x_0,x_1\},x_2,\{x_3,x_4\}],
\; [x_0, \{x_1,x_2\},\{x_3,x_4\}],\;\\ & &
[x_0\},\{x_1,x_2\},x_3,\{x_4 ],\;
 [x_0\},x_1,\{x_2,x_3\},\{x_4],
\end{eqnarray*}
where by, say $[\{x_0,x_1\},\{x_2,x_3\},x_4]$ we mean the point
$[x_0,x_1,x_2,x_3,x_4] \in \overline{M}_{0,5}$ with $x_0=x_1$ and
$x_2=x_3,$ and we let $t_{01,23}$ denote a tangent vector at that
point that maps to the previously defined tangent vectors on $X$
under the map that sends
$$[x_0,x_1,x_2,x_3,x_4] \to [x_0,x_2,x_3,x_4]$$
and the map that sends
$$
[x_0,x_1,x_2,x_3,x_4] \to [x_0,x_1,x_2,x_4].
$$
There are four different choices, however in the crystalline
setting the choice between these four points will not be
important; see the lemma below. Similarly we choose tangent
vectors at the four remaining basepoints with the same property.

{\bf 7.4. Frobenius.} Let $Y/\mathbb{Q}_{p}$ be a smooth variety,
and assume that there is a proper, smooth model
$\overline{\mathfrak{Y}}/\mathbb{Z}_p$   and a simple relative
normal crossings divisor $\mathfrak{D} \subseteq
\overline{\mathfrak{Y}}$ whose irreducible components are defined
over $\mathbb{Z}_p,$ with a fixed isomorphism
$\mathfrak{Y}_{\mathbb{Q}_p}\simeq Y,$ where
$\mathfrak{Y}:=\overline{\mathfrak{Y}}\setminus \mathfrak{D}.$
Using the isomorphism ${\rm Mic}_{uni}(Y/\mathbb{Q}_p)\simeq {\rm
Isoc}_{uni} ^{\dagger}(\mathfrak{Y}\otimes \mathbb{F}_p
/\mathbb{Z}_p) $ we obtain a frobenius action $F^*: {\rm
Mic}_{uni}(Y/\mathbb{Q}_p) \to {\rm Mic}_{uni}(Y/\mathbb{Q}_p)$
defined in [De] (see also [CS] and [\"{U}n]). And choosing
(tangential) basepoints $x$ and  $y$ with finite reduction  we
obtain a map $F_{*}: \, _{y} \mathcal{G}_{dR,x}(Y/\mathbb{Q}_p)
\to \, _{y}\mathcal{G}_{dR,x}(Y/\mathbb{Q}_p),$ where by
$\mathcal{G}_{dR}$ we denote the de Rham fundamental groupoid. In
fact the frobenius is independent of the choice of the model, but
we will not need this  below.

{\it p-adic integration.} Let $\mathbb{Q}_{p,st}$ denote the ring
of polynomials $\mathbb{Q}_{p}[l(p)],$ where $l(p)$ is a formal
variable that could be thought of as (a multi-valued) $\log p.$
Then $\log z:D(1,1^-) \to \mathbb{Q}_{p}$ uniquely extends to an
additive map  $\log z: \mathbb{Q}^{*} _{p} \to \mathbb{Q}_{p,st}$
such that  $\log p=l(p).$

Let $Y/\mathbb{Q}_p$ be a variety with a model
$\overline{\mathfrak{Y}}/\mathbb{Z}_p$ etc. as above and $x,y \in
X(\mathbb{Q}_p).$ Then Vologodsky [Vo] extending the work of
Coleman, Colmez, Besser, shows that there is  a canonical path
$_{y}c_{x} \in \, _{y}
\mathcal{G}_{dR,x}(Y/\mathbb{Q}_p)(\mathbb{Q}_{p,st})$ such that
$F_{*} (_{y}c_{x})=\,_{y}c_{x},$ and $_{y} c_{x}$ has image
$_{y}1_{x}$ under the canonical projection of
$_{y}\mathcal{G}_{dR,x}(\mathbb{Q}_{p,st})$ to
$\mathbb{Q}_{p,st}.$ $F$ coincides with the above when
 $x$ and $y$ have finite reduction with respect to the given model; in
fact  in this case  the element $_{y} c_{x}$ is defined over
$\mathbb{Q}_{p} \subseteq \mathbb{Q}_{p,st}.$

Similarly there is a canonical path $_{b}c_{a}$ even when $a$ and
$b$ are tangential basepoints. The same proof as in loc. cit. p.
17 extends to this case to show the existence and uniqueness of
the path satisfying the properties above.

 In order to describe this path, without
loss of generality  by the compatibility with respect to
concatenation, we will describe $_{u}c_{x}$ where  $x$ is a
genuine point and $u \in N^* _{D_y/\overline{Y}}(\mathbb{Q}_p) ,$
which is  not necessarily of finite reduction with respect to the
given model.

Let $\overline{Y}_{an}$ be the associated analytic space,
$\mathcal{U} \subseteq N_{D_{y}/\overline{Y}}$ a polydisc around
zero and $\varphi: \mathcal{U}_{log} \to \overline{Y}_{log,an}$ be
a map such that $\varphi(0)=y;$ the map on the monoids is the
identity map; the map induced by the differential at zero is the
identity map; and $\varphi$  is an closed immersion. Assume
without loss of generality that $x$ is in the image of $\varphi$
and $u$ is in $\mathcal{U}.$ Then we claim that the canonical
crystalline path between $\omega(x)\otimes\mathbb{Q}_{p,st} $ and
$\omega(u)\otimes \mathbb{Q}_{p,st}$ is given by $\lim_{\epsilon
\to 0} (\,_{u}c_{\epsilon}\cdot \,_{\varphi(\epsilon)}c_{x})$
where $_{u}c_{\epsilon }$ is  a path on
$N_{D_{y}/\overline{Y}}\setminus D_{(y)} ,$
$\,_{\varphi(\epsilon)}c_{x}$ is a path on $Y$ and
$\omega(\epsilon)$ and $\omega(\varphi(\epsilon))$ are identified
by the fact that the restriction to $\mathcal{U}$ of the pull-back
to the tangent space of a unipotent connection is its pull-back
via $\varphi$ as described above.

First in order to see that the above limit exists by choosing
coordinates and choosing a local trivialization
$(\mathcal{O}^{m}, d-\sum_{1\leq i \leq r}N_{i}\frac{dz_i}{z_i} )$
where $N_{i}$ are nilpotent, of the the underlying bundle of a
unipotent vector bundle with connection, we are reduced to showing
that: if
$$\varphi:D(0,1^{-})_{log} ^r \to D(0,1^-)_{log} ^{n}
$$
is a closed immersion of logarithmic analytic spaces, where both
of the spaces are endowed with the log structure associated to
$z_1 \cdots z_r=0,$ such that $\varphi(0)=0;$ and
$\frac{\varphi(z_i)}{z_i} \in 1 + \mathfrak{m}_{\overline{Y},0}$
then
$$\lim_{\epsilon \to 0} \prod_i exp(N_{i}\log
\frac{u_i}{\epsilon_i}) \cdot \prod_i exp(N_{i}\log
\frac{\varphi(\epsilon)_i} {x_i})$$
 exists.  Here $\log$ denotes
the multivalued extension of the logarithm described above. Since
$$\lim_{\epsilon \to 0}\log
\frac{\varphi(\epsilon)_i}{\epsilon_{i}}=0$$ the above limit
exists and is equal to $\prod_{i}exp(N_i \log\frac{u_{i}}{x_{i}}
)$. The standard arguments as in the section on tangential
basepoints show that the definition does not depend on the choice
of $\varphi.$ Let $\overline{c}$ denote this path just described.
The above argument also shows that if $(E,\nabla)$ is a unipotent
log isocrystal on $\overline{Y}_{log}$ with a local trivialization
$(\mathcal{O}^m,d-\sum N_{i}\frac{dz_i}{z_i})$ and $\varphi$ is a
map with the properties as above  from an open disc in
$N_{D_{y}/\overline{Y}}$ to $\overline{Y}$ then
$$
\lim_{v \to 0} \prod_{i} (exp(-N_{i} \log \frac{u_i} {v_i})
\,_{u}\overline{c}_{\varphi(v)})=1.
$$

In order to see that $_{u}\overline{c}_{x}$ is the canonical
crystalline path we need to show the invariance under frobenius.
Since the path is invariant under frobenius between genuine
basepoints to prove the invariance in general it suffices to prove
this invariance for the limit. We choose local coordinates and a
local trivialization of $(E,\nabla)$ as above.

Let $\mathcal{F}$ be a local lifting of the frobenius on the
special fiber to a neighborhood of $y$ in $ \overline{Y}_{log}$
relative to the model $\overline{\mathfrak{Y}},$ that fixes $y.$
Then the principal part $P(\mathcal{F})$ of $\mathcal{F}$ defines
the corresponding lifting of frobenius to
$\overline{N}_{D_{y}/\overline{Y},\, log}.$ Note that
$$F_{*} (_{u}
\overline{c}_{\varphi(v)})=\,_{u}{\rm par}_{P(\mathcal{F})(u)}\;
\mathcal{F}_{*}(_{u} \overline{c}_{\varphi(v)})\;
\,_{\mathcal{F}(\varphi(v))}{\rm par}_{\varphi(v)},$$ where the
$par$ on the left and on the right denote the parallel transport
along the connection on the spaces
$\overline{N}_{D_{y}/\overline{Y},log}$ and on
$\overline{Y}_{log}$ respectively. Note that since we are only
interested  in $_{\mathcal{F}(\varphi(v))}{\rm par}_{\varphi(v)}$
when $v$ tends to zero we may assume that the parallel transport
is in fact the parallel transport along the trivialized connection
on $\varphi(\mathcal{U})_{log}.$ In other words in the limit it
can be replaced with $\prod_{i}exp(N_{i} \log
\frac{\mathcal{F}(\varphi(v))_i}{\varphi(v)_i}).$

First note that since pulling back by $\mathcal{F}$ multiplies the
residues of a connection with  $p,$ by the above limit computation
for $_{u}\overline{c}_{\varphi(v)}$ we have
$$
\lim_{v \to 0} \prod_{i} exp(-pN_{i} \log \frac{u_i} {v_i})
\,\mathcal{F}_*( _{u}\overline{c}_{\varphi(v)})=1.
$$
Therefore
$$
\lim_{v \to 0} \prod_{i} exp\big(-N_{i} \log \frac{u_i ^p \,u_i\,
\mathcal{F}(\varphi(v))_i} {v_i ^p \, P(\mathcal{F})(u)_i \,
\varphi(v)_i}\big)\, F_*( _{u}\overline{c}_{\varphi(v)})=1.
$$
Note that $P(\mathcal{F})$ is given by
$P(\mathcal{F})(z_1,\cdots,z_{r} )=(a_{1}z_1 ^p, \cdots, a_r z_r
^p )$ for some $a_i \in 1+ p \mathbb{Z}_p.$ Since $\varphi$ is
identity on the tangent space $\lim_{v \to 0}
\frac{\mathcal{F}(\varphi(v))_i}{v_i ^p}=a_i=\lim_{v \to 0}
\frac{P(\mathcal{F})(u)_i}{u_i ^p},$ and $\lim_{v \to 0}
\frac{v_i}{\varphi(v)_i}=1 .$ And hence
$$
\lim_{v \to 0} \prod_{i} exp\big(-N_{i} \log \frac{u_i}
{v_i}\big)\, F_*( _{u}\overline{c}_{\varphi(v)})=1,
$$
and therefore
$F_{*}(_{u}\overline{c}_{\varphi(v)})=_{u}\overline{c}_{\varphi(v)}.$

This gives the description of the canonical crystalline path
between possibly tangential basepoints. If $x$ and $y$ have good
reduction relative to a model then $_{y} c_{x} \in \,
_{y}\mathcal{G}_{dR,x}(\mathbb{Q}_{p,st})$ is in fact defined over
$\mathbb{Q}_{p}.$ Again the argument cited above works in this
case if we note that  in the good reduction case in order to
define the frobenius, by the above method of tangential
basepoints, one does not have to tensor with $\mathbb{Q}_{p,st}.$

In order to see that changing the tangential basepoints by an
$r-$tuple of roots of unity will have no effect in the crystalline
de Rham theory we need the following lemma.

\begin{lemma}
Let $\mathfrak{Y}:=\mathbb{G}_{m} ^{r} /\mathbb{Z}_p ,$
$\overline{\mathfrak{Y}}$ be the standard compactification and $x$
and $y$ be in $\mathfrak{Y}(\mathbb{Z}_p).$ Let $$\{e_{i}\}_{1\leq
i \leq r }\in \pi_{1,dR}(Y,\omega(dR))$$ denote the residues at 0.
Then
$_{\omega(dR)}e(dR)_{y}F_{*}(_{y}e(dR)_{x})\,_{x}e(dR)_{\omega(dR)}=
\prod_{i} exp(e_i \log \frac{y_i ^{1-p}}{x_i ^{1-p}}).$
 \end{lemma}

{\it Proof.} We will assume without loss of generality that $x=1$
as the general statement above follows from this by concatenation
of paths. Let $\mathcal{F}(z_1,\cdots,z_r)=(z_1 ^{p},\cdots,z_r
^p)$ this is a lifting to $\overline{Y}_{log}$ of the frobenius on
the special fiber. Fix the beginning point as $1$ and consider the
$\pi_{1,dR}(Y,1)-$torsor of paths $_{y}\mathcal{G}_{1,dR}$ as $y$
varies. Note that $\mathcal{F}_{*}$ defines a horizontal map
$\mathcal{F}_{*}: _{y}\mathcal{G}_{1,dR} \to \mathcal{F}^{*}\,
_{y}\mathcal{G}_{1,dR},$ where $\mathcal{G}_{dR}$ is endowed with
its canonical connection. This gives a differential equation for
$\mathcal{F}_{*}(_{y}e(dR)_{1})$ and solving this we find
$\mathcal{F}_{*}(_{y}e(dR)_{1})=1,$ see [\"{U}n] for a similar
computation where more details are given.  Since
$F_{*}(_{y}e(dR)_{1})=\,_{y}{\rm par}_{\mathcal{F}(y)}
\mathcal{F}_{*}(_{y}e(dR)_{1}),$ as $\mathcal{F}(1)=1$ and
$_{y}{\rm par}_{y^p}=\prod_{i} exp(e_{i}\log y^{1-p} )$ the
statement follows. Note that if $y_{i}/x_{i}$ are roots of unity
for all $i$  then the above shows that
$F_{*}(_{y}e(dR)_{x})=\,_{y}e(dR)_{x},$ as $\log(a)=0$ if $a$ is a
root of unity.  \hfill $\Box$

{\bf 7.5. p-adic multi-zeta values.}

{\it Notation.} For a smooth $Y/K,$ and $y\in Y(K)$ let
$\mathcal{U}_{dR}(Y,y)$ denote the universal enveloping algebra of
${\rm Lie}\pi_{1,dR}(Y,y)$ and $\hat{\mathcal{U}}_{dR}(Y,y)$ be
its completion with respect to its augmentation ideal. It is a
cocommmutative Hopf algebra  and its topological dual is the Hopf
algebra of functions on $\pi_{1,dR}(Y,y).$

Letting $e_0,e_1,$ and $e_{\infty} \in {\rm
Lie}\pi_{1,dR}(X/\mathbb{Q}_p,t_{01})$  denote the residues
corresponding to the points $0,1, \infty$ in $\overline{X}$
respectively, ${\rm Lie}\pi_{1,dR}(X,t_{01})\simeq {\rm Lie}
<<e_0,e_1 >>.$  $\hat{\mathcal{U}}_{dR}(X,t_{01})$ is isomorphic
to the ring of associative formal power series on $e_0$ and
$e_{1}$ with the coproduct $\Delta$ given by $\Delta(e_0)=1\otimes
e_0+e_{0}\otimes 1, $  and $\Delta(e_1)=1 \otimes e_1 + e_1
\otimes 1. $ By the duality above $\mathbb{Q}_p -$rational points
of $\pi_{1,dR}(X,t_{01})$ correspond to associative formal power
series  $a$ in $e_0$ and $e_1$ with coefficients in
$\mathbb{Q}_p,$ that start with 1 and satisfy $\Delta(a)=a\otimes
a.$

We let $$g:=_{t_{01}}e(dR)_{t_{10}} F_{*}(_{t_{10}}e(dR)_{t_{01}})
\in \pi_{1,dR}(X,t_{01}).$$ This is the series that defines the
p-adic multi-zeta values. In particular the coefficient of the
term $e_0 ^{s_{k}-1}e_1\cdots e_0 ^{s_1 -1}e_1$ is by definition
$p^{\sum s_i}\zeta_{p}(s_k,\cdots,s_1).$ These values also
determine $g.$

{\bf 2-cycle relation.}  Let $\gamma:=_{t_{10}}e(dR)_{t_{01}},$
and we put $g:=g(e_0,e_1)=\gamma^{-1}F_{*}(\gamma) \in
\pi_{1,dR}(X_{\mathbb{Q}_{p}},t_{01}).$  We would like to see that
$$
g(e_1,e_0)g(e_0,e_1)=1.
$$
Let $\tau$ be the automorphism of $X$ that maps $z$ to $1-z.$ Then
$$
\tau(t_{01})=t_{10}, \, \tau_{*}(\gamma)=\gamma^{-1}, \,
\tau_{*}(e_0)=\gamma e_1\gamma^{-1}, \, \tau_{*}(e_1)=\gamma e_0
\gamma^{-1}.
$$
We have
$$
\gamma g(e_1,e_0)
\gamma^{-1}=\tau_{*}(g(e_0,e_1))=\tau_{*}(\gamma^{-1})\tau_{*}(F_{*}\gamma)=
\tau_{*}(\gamma^{-1})F_{*}(\tau_{*}(\gamma))=\gamma
F_{*}(\gamma^{-1}).
$$
Therefore
$$
g(e_1,e_0)=F_{*}(\gamma^{-1})\gamma=g(e_0,e_1)^{-1}.
$$

{\bf 3-cycle relation.} Let $\delta:=_{t_{\infty
0}}e(dR)_{t_{01}},$ $r:=_{t_{1 \infty }}e(dR)_{t_{10}},$ and
$q:=r^{-1} \gamma=_{t_{1 \infty}}e(dR)_{t_{01}}.$ And let
$e_{\infty}:=\delta^{-1}{\rm res}_{t_{\infty 0}}\delta $ be the
Lie element describing the residue at $\infty$ with basepoint
$t_{01}. $ Then we would like to see that
$$
g(e_{\infty},e_0)g(e_1,e_{\infty})g(e_0,e_1)=1.
$$
Let $\omega$ be the automorphism of $X$ that sends $z$ to
$\frac{1}{1-z}.$ Then
\begin{eqnarray*}
\omega(t_{01})&=&t_{1 \infty}, \,{\rm and } \, \,
\delta=\omega_{*}(q)q,\, \omega_*(e_0)=q e_1 q^{-1},\,
\omega_*(e_1)=qe_{\infty}q^{-1},\\
\omega^2_*(e_0)&=&\omega^2_*(q)^{-1}e_{\infty}\omega^2_*(q), \,
\omega^2_*(e_1)=\omega^2_*(q)^{-1}e_0\omega^{2}_*(q)
\end{eqnarray*}
Applying frobenius to
$$
\omega^{2}_*(q)\omega_*(q)q=1
$$
we obtain
\begin{eqnarray*}
1&=&
\omega^{2}_*(F_{*}q)\omega_*(F_{*}q)F_{*}q=\omega^{2}_*(F_{*}q)\omega_*(F_{*}q)
q q^{-1}F_{*}q\\
&=&\omega^{2}_*(F_{*}q)\omega_*(q)q.\,
q^{-1}\omega_*(q)^{-1}\omega_*(F_{*}q)q.\, g\\
&=&\omega^{2}_*(F_{*}q)\omega_*(q)q.\,
q^{-1}\omega_{*}(g(e_0,e_1)) q.\,
g=\omega^{2}_*(F_{*}q)\omega_*(q)q. \, g(e_1,e_{\infty}).\,
g(e_0,e_1)\\
&=&\omega^2_{*}(q) \omega^{2}_*(g(e_0,e_1))\omega^2_{*}(q)^{-1}.
\, \omega^{2}_*(q)\omega_*(q)q. \, g(e_1,e_{\infty}).\,
g(e_0,e_1)\\
&=& g(e_{\infty},e_0)g(e_1,e_{\infty})g(e_0,e_1)
\end{eqnarray*}
since $r^{-1}F_{*}(r)=1,$ since $\log (-1)^{1-p}=0.$

{\bf 5-cycle relation.} In this section we identify $M_{0,5}$ with
$X^{2}\setminus \{(z,w)|zw=1 \} $ by the map $(z,w) \to
[0,z,1,w^{-1},\infty].$ By this identification let
$t_{ij}(i,\epsilon)$ denote the tangent vector $(t_{ij},0)$ at the
point $(i,\epsilon).$ And we choose, $t_{01,34}$ to be the tangent
vector $(1,1)$ at the point $(0,0)$ and $t_{12,34}$ to be the
vector $(-1,1)$ at the point $(1,0).$

\begin{lemma}
We have
$$
\lim_{y \to 0} exp(-e_{34} \log y  )\,
_{t_{01}(0,y)}c_{t_{01,34}}=1,
$$
where  we use the canonical de Rham paths to identify the
different basepoints.
\end{lemma}

{\it Proof.} By the description of the crystalline path with
endpoints at a tangent vector we have, after always using the de
Rham trivialization to identify the different basepoints,
$$
\lim_{(x,y)\to (0,0)} exp(-e_{34}\log y )exp(- e_{01} \log x )
\,_{(x,y)} c_{t_{01,34}}=1
$$
and
$$
\lim_{x \to 0} exp(-e_{01}\log x )\,_{(x,y)} c_{t_{01}(0,y)}=1,
$$
which give the statement in the lemma. \hfill $\Box$

Similarly,
$$
\lim_{y \to 0} exp(-e_{34} \log y  )\,
_{t_{10}(1,y)}c_{t_{12,34}}=1.
$$

Therefore we have
$$
_{t_{12,34}}c_{t_{01,34}}=\lim_{\epsilon \to 0} exp(-e_{34}\log
\epsilon)\cdot \; _{t_{10}(1,\epsilon)}c_{t_{01}(0,\epsilon)}\cdot
\; exp(e_{34}\log \epsilon).
$$

Because of good reduction we know that the left hand side in fact
is defined over $\mathbb{Q}_{p}$ therefore it is unchanged if we
basechange by the map $\mathbb{Q}_{p,st} \to \mathbb{Q}_{p}$ that
sends $l(p)$ to 0. We do this in the remaining part of the
section. Therefore restricting $\epsilon$ to powers of $p$ we
obtain
$$
_{t_{12,34}}c_{t_{01,34}}=\lim_{N \to \infty} \,
_{t_{10}(1,p^N)}c_{t_{01}(0,p^N)}.
 $$

With the coordinates as above let $X_{\epsilon}$   be the
subvariety of $M_{0,5}$ defined by $w-\epsilon=0,$ for $\epsilon
\in \mathbb{Q}_{p} ^{*}.$ Therefore $X_{\epsilon} \simeq X
\setminus \{1/\epsilon \}.$ Note that
$_{t_{01}(1,\epsilon)}c_{t_{01}(0,\epsilon)}$ is the image  of the
analogous element  $$c_{\epsilon}:=_{10}c_{t_{01}}(X_{\epsilon})$$
under the inclusion $X_{\epsilon} \to M_{0,5}.$

We will need the following
\begin{lemma}
Let $a$ be a monomial in $e_{0},$ $e_1$ and $e_{p^{-N}}$ such that
$a$ contains an $e_{p^{-N}}.$ If we denote the image of $c_{p^N}$
under the canonical maps
$$
_{t_{10}}\mathcal{G}_{t_{01},dR}(X_{p^N}(\mathbb{Q}_{p,st}) \simeq
\pi_{1,dR}(X_{p^N},\omega(dR)) \subseteq
\mathcal{U}_{dR}(e_0,e_1,e_{p^{-N}})(\mathbb{Q}_{p,st})\simeq
\mathbb{Q}_{p,st}<<e_0,e_1,e_{p^{-N}} >>,
$$
by the same symbol, then $\lim_{N \to \infty} c_{p^N}[a]=0.$
\end{lemma}

{\it Proof.} First note that  since $c_{p^N}$ is a group-like
element of $\hat{\mathcal{U}}_{dR}$ the coefficient $c_{p^N}[a]$
is a (finite) linear combination of terms $c_{p^N}[b]$ where $b$
is a monomial that   ends with $e_{p^{-N}}.$ In order to see this
first we write $a:=a' \cdot a_{0},$ where $a'$ ends with
$e_{p^{-N}}$ and $a_{0}$ does not contain any $e_{p^{-N}}.$ Then
we  compare the coefficients of the term $a' \otimes a_{0}$ on
both sides of the equality $\Delta(c_{p^N})=c_{p^N}\otimes
c_{p^N},$ and use induction on the number of terms on the right of
the last $e_{p^{-N}}$ in $a.$

Therefore  without loss of generality we assume that  $a$  ends
with $e_{p^{-N}}.$ For $i \in \{0,1, p^{-N} \},$ let
$\omega_{i}:=d \log (z-i).$ For a sequence of $\omega_i's:$
$\omega_{i_r}, \dots, \omega_{i_1}$ with $i_1\neq 0,$ we associate
the iterated integral $ I_{i_r,\cdots,i_1}(z_r):= \int_{0}
^{z_{r}} \omega_{i_r} \circ \cdots \circ \omega_{i_1} $ defined
successively by
$$I_{i_r,\cdots,i_1}(z_r):=\int_{0} ^{z_r}
\omega_{i_r}(z_{r-1})I_{i_{r-1},\cdots ,i_1}(z_{r-1}),$$ where by
the misleading notation $\int_{0} ^{z}$ we mean the
anti-derivative of the integrand in the sense of the Coleman
integral starting at the point 0. In the Coleman integration  we
let $\log p=l(p)=0,$ to be compatible with the above choice of a
branch. This iterated integral is a locally analytic function on
$X_{p^{N},an}.$ Note that since we assume that $\omega_{i_1}$ does
not have a singularity at zero the iterated integral is in fact
analytic around 0 with value 0 at 0. By definition of p-adic
integration the iterated integral $I_{i_r,\cdots,i_1}(z)$ is the
coefficient $_{z}c(X_{p^N})_{t_{01}}[e_{i_r} \cdots e_{i_1}].$ By
the construction of Coleman integration
$$I_{i_r,\cdots,i_1}(z)=\sum_{0 \leq k \leq n } f_{k}(z) \log ^k
(z-1) $$ on  $U\setminus \{ 1\},$ where $U$  is a neighborhood of
1 and $f_{k}$ are analytic on $U$ [Co].

By the definition of the tangential basepoints and the description
above
$$_{t_{10}}c_{t_{01}}[a]=\lim_{z \to 1}\Big(  exp(-e_{1} \log (z-1)) \cdot \,
 _{z}c_{t_{01}}\Big) [a].$$
Therefore we obtain
$$
 _{t_{10}}c_{t_{01}}[a]=\lim_{M \to \infty} f_{0}(1+p^M).
$$

Therefore in order to prove the lemma it suffices to show that if
$(\varphi_{N}(z))_{N}$ is a sequence of analytic functions on some
$D(0,(1+\delta)^{-})$ which converge uniformly to 0 with
$\varphi_{N}(0)=0,$ then we have
$$
\lim_{N_{\to \infty}} \lim_{M \to \infty} \int_{0} ^{1+p^M}
\omega_{i_k} \circ \cdots \circ \omega_{i_1} \circ d
\varphi_{N}(z)=0.
 $$

We show this by induction on the weight, i.e. the number of terms
in the iterated integral. The assertion is clear if the weight is
one. Assume that the weight is greater than one. If
$\omega_{i_1}=d\log(z-p^{-N})$ then by noting that $$\int_{0} ^{z}
d\log(z-p^{-N})\circ d\varphi_{N}(z)$$ satisfies the conditions
for $\varphi_{N},$ we reduce to the case with one lower weight.
The same is true if $\omega_{i_1}=d\log z.$ So assume that
$\omega_{i_1}=d \log (z-1).$ Note that $$\int_{0} ^{z} d \log
(z-1) \circ d \varphi_{N}(z)=\int_{0} ^{z} d \log(z-1)\cdot
(\varphi_{N}(z)-\varphi_{N}(1) ) + \varphi_{N}(1)\cdot \log
(z-1)$$

The first term on the right satisfies the same conditions as
$\varphi_{N}$ therefore we only need to take care of the term
$\varphi_{N}(1)\log(z-1).$ In other words need to show that
$$
\lim_{N \to \infty} \lim_{M \to \infty} \varphi_{N}(1) \cdot
\int_{0} ^{1+p^M} \omega_{i_k} \circ \cdots \circ \omega_{i_1}=0.
$$
If there is a $d \log (z-p^{-N})$ among the $\omega_{i_j}$ this
follows from the induction hypothesis, otherwise the iterated
integral does not depend on $N$ and  $\lim_{N \to
\infty}\varphi_{N}(1)=0$ implies the assertion. \hfill $\Box$

This lemma implies that $_{t_{12,34}}c_{t_{01,34}}$ only consists
of the terms $e_{01}$ and $e_{12}.$  Let $$\alpha:=\, _{t_{01,34}}
e(dR) _{t_{12,34}} \cdot \, _{t_{12,34}}c_{t_{01,34}}.$$ Then
$F_{*}\alpha=( _{t_{12,34}}g_{t_{01,34}})^{-1} \cdot
\,_{t_{12,34}}c_{t_{01,34}} $ and $F_{*}(\alpha)$  only consists
of the terms $e_{01}$ and $$(_{t_{12,34}}g_{t_{01,34}})^{-1}\cdot
e_{12}\cdot \, _{t_{12,34}}g_{t_{01,34}}.$$ From this, by
induction on the weight, we see that $_{t_{12,34}}g_{t_{01,34}}$
also consists only of $e_{01}$ and $e_{12}.$ By functoriality
$_{t_{12,34}}g_{t_{01,34}}$ maps to $_{t_{10}}g_{t_{01}}$ under
the map $M_{0,5} \to X$ and hence we have
$_{t_{12,34}}g_{t_{01,34}}=g(e_{01},e_{12}).$ Similarly we get
\begin{eqnarray*}
 _{t_{12,13}}g_{t_{12,34}}&=&g(e_{34},e_{40})\\
_{t_{13,23}}g_{t_{12,13}}&=&g(e_{12},e_{23})\\
 _{t_{01,23}}g_{t_{13,23}}&=&g(e_{40},e_{01})\\
_{t_{01,34}}g_{t_{01,23}}&=&g(e_{23},e_{34})\\
\end{eqnarray*}

Similar to the case of other  relations we find that
$$
_{t_{01,34}}g_{t_{01,23}} \cdot\,_{t_{01,23}}g_{t_{13,23}}\cdot\,
_{t_{13,23}}g_{t_{12,13}} \cdot\, _{t_{12,13}}g_{t_{12,34}}
\cdot\, _{t_{12,34}}g_{t_{01,34}}=1,
$$
which gives  the Drinfel'd-Ihara relation
$$
g(e_{23},e_{34})\cdot g(e_{40},e_{01})\cdot  g(e_{12},e_{23})
\cdot g(e_{34},e_{40})\cdot g(e_{01},e_{12})=1.
$$

\end{document}